# ENDOSCOPIC CHARACTER IDENTITIES FOR METAPLECTIC GROUPS

CAIHUA LUO

ABSTRACT. In this paper, we prove the conjectural endoscopic character identities for tempered representations of metaplectic group $Mp_{2n}$ based on the formalism of endoscopy theory by J. Adams, D. Renard and W.W. Li.

## 1. INTRODUCTION

Recently, in a series of papers [Li12b, Li12a, Li14b, Li13, Li14a], W.W. Li has established the invariant trace formula for finite central covering groups in the spirit of the formulation of Arthur's trace formula. In particular, motivated by J. Adams and D. Renard's work on developing an endoscopy theory for $Mp(2n, \mathbb{R})$, W.W. Li (see [Li11, Li15]) has built up an almost complete endoscopy theory which places $Mp(2n)$ in the framework of the Langlands program. On the other hand, Howe's theta correspondence theory plays an essential role in the study of (automorphic) representations of $Mp(2n)$ and has produced fruitful developments. More specifically, the local Langlands correspondence for $Mp(2n)(\mathbb{R})$ has been established by J. Adams and D. Barbasch in [AB98], and almost 20 years later, the $p$-adic version was obtained by Gan and Savin in [GS12]. In view of these, parallel to Arthur's work in his monumental book [Art13], it is not surprising to expect and explore the following two fundamental "theorems":

- Endoscopic character identities for $Mp(2n)$,
- Multiplicity formula of $L^2_{disc}(Sp(2n, \dot{F})\backslash Mp(2n, \mathbb{A}))$ (Arthur's conjecture).

Along the lines of Arthur's standard model argument, W.W. Li in [GL16] proposed a strategy to tackle the above two problems, and formulated a conjecture for the desired form of the endoscopic character identities which is compatible with Waldspurger's influential work on $Mp(2)$ (cf. [Wal80, Wal91]). On the other hand, based on the powerful tool of theta correspondence, Gan and Ichino has worked out the elliptic tempered part of $L^2_{disc}(Sp(2n)(\dot{F})\backslash Mp(2n)(\mathbb{A}))$ in the sense of Arthur's $A$-packet (cf. [GI17]), where $\dot{F}$ is a global field and $\mathbb{A} = \mathbb{A}_{\dot{F}}$ is the associated *Adèle* ring. Note that the endoscopic character identity over $\mathbb{R}$ has been proved by D. Renard (cf. [Ren99]), but it has not yet been expressed in the form of characters of component groups. Note also that the elliptic stable trace formula has been laid down by W.W. Li in [Li15] and the spherical fundamental lemma has been established in [Luoar]. In view of these facts, we may apply a local-global argument via Arthur's stable multiplicity formula to tackle the endoscopic character identity problem once the real case is reformulated to fit into the global argument. In fact, this is exactly what we will do in this paper. Herein we would like to mention that J. Schultz's PhD thesis has established the endoscopic character identities for $Mp(2)$ (cf. [Sch98]), and T. Howard has proved the compatibility property of parabolic induction for the principal endoscopic group (cf. [How10]). In what follows, let us briefly recall the conjectural endoscopic character identities for $Mp(2n)$ formulated by W.W. Li in [GL16].

Let $F$ be a local field of characteristic 0. Fix a non-trivial character $\psi$ of $F$ and let $\widetilde{Sp}^{(8)}(W)$ be Weil's metaplectic 8-fold covering of $Sp(W)$ which is a pushout of the metaplectic 2-fold cover

$Mp(W)$ via $\mu_2 \hookrightarrow \mu_8$, i.e.

$$\begin{CD} \mu_2 @>>> Mp(W) @>p>> Sp(W) \\ @VVV @VVV @| \\ \mu_8 @>>> \widetilde{Sp}^{(8)}(W) @>p>> Sp(W) \end{CD}$$

Let $G = Sp(W)$, $\tilde{G} = \widetilde{Sp}^{(8)}(W)$ with $\dim W = 2n$. As shown by J. Adams, D. Renard and W.W. Li, the elliptic endoscopic groups of $\tilde{G}$ are the split orthogonal groups associated to the ordered pairs $(n', n'')$ with $n' + n'' = n$:

$$H = H_{n',n''} = H' \times H'' := SO(2n'+1) \times SO(2n''+1).$$

On the other hand, theta correspondence gives rise to the following bijective maps obtained by Gan–Savin for non-archimedean fields and Adams–Barbasch for the real archimedean field, i.e.

$$\theta_\psi : Irr\ \tilde{G} \longleftrightarrow \bigsqcup_{V_n} Irr\ SO(V_n)$$

where

- $V_n$ runs over all $2n+1$-dimensional quadratic spaces of fixed discriminant over $F$;
- $Irr\ \tilde{G}$ is the set of isomorphism classes of irreducible genuine representations of $\tilde{G}$, and $Irr\ SO(V_n)$ is the set of isomorphism classes of irreducible representations of $SO(V_n)$.

Combining with the local Langlands correspondence (abbreviated as LLC) for $SO(2n+1)$ (due to Arthur and Moeglin for non-archimedean fields, and to Langlands and Shelstad for archimedean fields), one gets a LLC for $\tilde{G}$:

$$Irr\ \tilde{G} \longleftrightarrow \{(\phi, \eta)\}$$

where

- $\phi : WD_F \longrightarrow Sp_{2n}(\mathbb{C})$ is an $L$-parameter for $SO(2n+1)$;
- $\eta \in Irr\ A_\phi$ with $A_\phi = \pi_0(Z_{Sp_{2n}}(Im(\phi)))$ the component group of $\phi$.

Thus, given a tempered $\phi$, we have an associated local $L$-packet

$$\Pi_\phi^{\tilde{G}} = \{\tilde{\pi}_\eta : \eta \in Irr\ A_\phi\}.$$

Now we may state the conjectural endoscopic character identities as follows.

**Main Theorem.** *Given a tempered L-parameter $\phi$ of $\tilde{G}$ and an $s \in Z_{Sp(2n)}(Im(\phi))$ with $s^2 = 1$, so that $s$ determines an elliptic endoscopic group $H_s = SO(2n'+1) \times SO(2n''+1)$ together with an L-parameter $\phi_{H_s} = \phi' \times \phi''$ of $H_s$, one has an associated L-packet $\Pi_{\phi' \times \phi''}^{H_s}$ and may then form the stable distribution*

$$\Theta_{\phi' \times \phi''} = \sum_{\pi \in \Pi_{\phi' \times \phi''}^{H_s}} \Theta_\pi$$

*where $\Theta_\pi$ is the normalized character distribution of $\pi$. Then we have the following local character identity:*

$$Trans_{H_s}^{\tilde{G}}(\Theta_{\phi' \times \phi''}) = \epsilon(1/2, \phi'', \psi) \sum_{\eta \in \hat{A}_\phi} \eta(s) \Theta_{\tilde{\pi}_\eta},$$

*where $Trans_{H_s}^{\tilde{G}}$ is the endoscopic transfer of stable distribution of $H_s$ to genuine invariant distribution of $\tilde{G}$.*

Let us end the introduction by giving a brief outline of this article. In section 2, we will recall the local Langlands correspondence for $Mp(2n)$ via theta correspondence from $SO(2n+1)$, and check the compatibility of Shelstad and Arthur's local Langlands correspondence for discrete series representations of $SO(2n+1)(\mathbb{R})$ with J. Adams' version (cf [Ada10]). The LLC for $Mp(2n)$ plays an important role in the procedure of reinterpreting D. Renard's results in terms of characters of component groups. In section 3, we will recall the endoscopy theory laid down by W.W. Li, and the stable trace formula for $H_s$ (due to Arthur) and the stable elliptic trace formula for $\tilde{G}$ (due to W.W. Li) which we would use later on. The last three sections are devoted to the proof of the Main

Theorem.

**Acknowledgements.** We are much indebted to Professor Wee Teck Gan for his guidance and numerous discussions on various topics. We would like to thank Professors Jeffrey Adams and Wen-Wei Li for discussions on this paper during a conference at CIRM, Luminy, France. We would also like to thank Professor Atsushi Ichino for inviting me to visit Kyoto University to give a talk on this paper. Thanks are also due to the referee for his/her detailed comments.

## NOTATIONS AND FACTS

- Let $F$ be a local field of characteristic 0, and $\psi_F : F \to S^1$ a non-trivial character of conductor $\mathcal{O}_F$.
- Let $W_F$ be the Weil group of $F$, and $WD_F = W_F \times SL_2(\mathbb{C})$ the Weil-Deligne group.
- Let $\dot{F}$ be a number field, and $\mathbb{A}$ the associated *Adèle* ring. Fix a non-trivial character $\psi : \dot{F}\backslash\mathbb{A} \to S^1$, with $\psi = \otimes_v \psi_v$. Denote by $V_{\dot{F}}$ ($V_\infty$ *resp.*) the set of (archimedean resp.) places of $\dot{F}$.
- Let $W$ be a symplectic vector space and fix a complete polarization $W = X + Y$ with associated standard sympletic basis of $W$, we then obtain an identification

$$Sp(W) = \{g = \begin{pmatrix} a & b \\ c & d \end{pmatrix} \in GL_{2n}(F) \mid g \begin{pmatrix} & 1 \\ -1 & \end{pmatrix} {}^t g = \begin{pmatrix} & 1 \\ -1 & \end{pmatrix}\}.$$

- To reserve $Mp(W)$ for Weil's twofold covering group of $G = Sp(W)$, we denote $\tilde{G} = Mp(W) \times_{\mu_2} \mu_8$ to be Weil's eightfold covering group of $Sp(W)$. Conventionally, for $\tilde{\delta} \in \tilde{G}$, we always write $\delta$ for $p(\tilde{\delta})$. Regarding $\tilde{G} = G \times \mu_8$ as a set given by the Rao cocycle (cf. [Kud96, Theorem 4.5]), we write $\tilde{\delta} = (\delta, \epsilon(\tilde{\delta}))$, with $\epsilon(\tilde{\delta}) \in \mu_8$.
- For the hyperspecial compact subgroup $K = Sp(2n)(\mathcal{O}_F)$ of $Sp(2n)$, $K$ splits in $\tilde{G}$ if the residue characteristic of $F$ is not 2. Globally, $Sp(2n, \dot{F})$ splits in $\tilde{G}(\mathbb{A})$.
- For $\tilde{x}, \tilde{y} \in \tilde{G}$, $\tilde{x}$ and $\tilde{y}$ commute if and only if $x$ and $y$ commute in $G$.
- For a maximal split torus $T$ of $G$, let $\tilde{T}$ be the preimage of the projection map $p$, one may then define a genuine $W^G(T)$-invariant character $\chi_{\psi_F}$ of $\tilde{T}$ as in [Kud96] as follows.

$$\chi_{\psi_F} : \left(\begin{pmatrix} a & \\ & {}^t a^{-1} \end{pmatrix}, \epsilon\right) \longmapsto \epsilon \gamma(det(a), \psi_F)^{-1}$$

where $\gamma(\cdot, \psi_F)$ is the relative Weil index. Note that this is also compatible with the natural splitting of $T$ in $\tilde{T}$ given by $\sigma_Y$ in [Kud96, Theorem 4.5]

$$\sigma_Y : g = \begin{pmatrix} a & \\ & {}^t a^{-1} \end{pmatrix} \longmapsto (g, \gamma(det(a), \psi_F)),$$

i.e. $\chi_{\psi_F} \circ \sigma_Y = id$. Throughout the paper, we will write $det(\gamma)$ in place of $det(a)$ with $\gamma = diag\{a, {}^t a^{-1}\}$ for simplicity.
- Let $Irr\,\tilde{G}$ be the set of isomorphism classes of genuine irreducible admissible representations of $\tilde{G}$, and $Irr\,H$ the set of isomorphism classes of irreducible admissible representations of $H$.
- For Weyl group, we sometimes abbreviate $W_\mathbb{C}^G$ for $W(G(\mathbb{C}), T(\mathbb{C})) = N_{G(\mathbb{C})}(T(\mathbb{C}))/T(\mathbb{C})$, and $W_\mathbb{R}^G$ for $W(G(\mathbb{R}), T(\mathbb{R})) = N_{G(\mathbb{R})}(T(\mathbb{R}))/T(\mathbb{R})$ when there is no danger of confusion.

## 2. LOCAL LANGLANDS CORRESPONDENCE AND MULTIPLICITY FORMULA

In this section, we first briefly recall the local Langlands correspondence (abbreviated as LLC) for $SO(2n+1)$ in the sense of Vogan packets, and the LLC for $\tilde{G}$ via theta correspondence. Then, we would like to check that Shelstad and Arthur's LLC for real forms of $SO(2n+1)$ are the same as Adams' version. Given the LLC, we will at last recall Gan–Ichino's global multiplicity formula for elliptic tempered $A$-parameters of $\tilde{G}$ in $L^2_{disc}(Sp(W)\backslash\tilde{G}(\mathbb{A}))$ which will play an essential role in our proof of the Main Theorem later on.

2.1. **LLC for $\tilde{G}$: Coarse version.** Let $H = SO(V) = SO(2n+1)$ be split, so that $\hat{H} = Sp(2n, \mathbb{C})$ is endowed with trivial $W_F$-action. Note that a key feature of Vogan's LLC is to treat representations $\pi$ of all pure inner forms $H'$ of $H$ simultaneously. Recall that the pure inner forms of $H$ are the groups $H'$ over $F$ which are obtained from $H$ via inner twisting by elements in the finite pointed set $H^1(F, H)$. As is well-known, the pure inner forms of $H$ are the groups $H_m = SO(V_m)$, where $V_m$ is an orthogonal space over $F$ with $dim(V_m) = 2n+1$ and the same discriminant as $V$. More precisely,

- If $F$ is non-archimedean and $n \geq 1$, there is a unique non-split pure inner form $H'$. In this case, we denote the split $H$ as $SO(V^+)$, while the non-split $H'$ as $SO(V^-)$.
- If $F = \mathbb{R}$ and $H = SO(p, q)$, then the pure inner forms are the groups $H' = SO(p', q')$ with $q \equiv q' \pmod{2}$ and $p + q = p' + q' = 2n + 1$.

An $L$-parameter for $H$ is an admissible homomorphism ([Bor79])

$$\phi : WD_F \longrightarrow \hat{H}.$$

We say $\phi$ is tempered if $\phi(W_F)$ is bounded. Further, if $\phi$ does not factor through any Levi subgroup of $\hat{H}$, we say $\phi$ is discrete (or elliptic). Such parameters are what we will focus on in this article. Given $\phi$, we define the S-group $S_\phi$ and the corresponding component group $A_\phi$ as follows.

$$S_\phi = Z_{\hat{H}}(Im(\phi)), \qquad A_\phi = \pi_0(S_\phi).$$

Also define the associated $L$-packet as $\Pi_\phi^H$. Now we may state Vogan's LLC as follows.

**LLC for $SO(2n+1)$.** *There exists a bijective correspondence*

$$L_V : \bigsqcup_{V_m} Irr\ SO(V_m) \longleftrightarrow \bigsqcup_{\phi \in \Phi} \{(\phi, \eta) : \eta \in Irr\ A_\phi\}$$

*where*

- *$V_m$ runs over all $2n+1$-dimensional quadratic spaces over $F$ of fixed discriminant;*
- *$\Phi$ is the set of equivalence classes of admissible representations of $WD_F$ in $Sp(2n, \mathbb{C})$.*

*Moreover, the correspondence $L_V$ preserves the property of temperedness and discreteness.*

*To be precise, suppose a tempered $\phi = \phi_{GL} \oplus \phi_0 \oplus \phi_{GL}^\vee$ factors through a Levi subgroup $\hat{M} = \prod_{i \in I} GL(n_i, \mathbb{C}) \times Sp(2n_0, \mathbb{C})$ of the dual group $\hat{H}$, where $\phi_{GL}$ and $\phi_0$ are the discrete parameters of $\prod_{i \in I} GL(n_i)$ and $SO(2n_0 + 1)$ respectively. Denote $M = \prod_{i \in I} GL(n_i) \times SO(2n_0 + 1)$, and choose a parabolic subgroup $P = MN$ of $H$ with Levi factor $M$, then*

$$\Pi_\phi^H = \bigsqcup_{\sigma \in \Pi_\phi^M} \{\text{the irreducible constituents of } I_P(\sigma)\}.$$

*Furthermore, for $\sigma_1, \sigma_2 \in \Pi_\phi^M$, the associated R-groups $R_{\sigma_1}$ and $R_{\sigma_2}$ of normalized induced representations from $\sigma_1$ and $\sigma_2$ respectively have the following good properties:*

- *$R_{\sigma_1}$ and $R_{\sigma_2}$ are abelian groups;*
- *$R_{\sigma_1} \simeq R_{\sigma_2}$.*

*In particular, $I_P(\sigma)$ is multiplicity-free.*

*Proof.* Note that this has been shown in Arthur's monumental book (see [Art13]) for $p$-adic quasi-split $SO(V_m)$, while the real case has been established by Adams–Barbasch–Vogan in [ABV12]. Recently the non-split case has also been established in [MR17]. As for R-group, one may refer to [Art13][Chapters 6.5 & 6.6] and [CG16][Theorem 3.9] for more details. $\square$

Note that J. Adams has established an explicit labeling on the discrete part of LLC for $SO(2n+1)$. For our purpose, we need to check the compatibility of these two versions which will be carried out in the next Section 2.2.

On the other hand, Howe's theory of theta correspondence theory provides us a bijection between $Irr\ \tilde{G}$ and $\bigsqcup_{V_m} Irr\ SO(V_m)$, which was established by Adams–Barbasch for real case (cf. [AB98]), and later proved by Gan–Savin for $p$-adic case (cf. [GS12]). Building upon this, we could obtain a LLC for $\tilde{G}$ as follows.

**LLC for $\tilde{G}$.** *Given the LLC for $SO(2n+1)$ in the sense of Vogan, taking $\hat{\tilde{G}} = \hat{H}$, then there is a bijective correspondence*

$$L_\Theta : Irr\ \tilde{G} \overset{\Theta_\psi}{\longleftrightarrow} \bigsqcup_{V_m} Irr\ SO(V_m) \overset{L_V}{\longleftrightarrow} \bigsqcup_{\phi \in \Phi} \{(\phi, \eta) : \eta \in Irr\ A_\phi\}$$

*where the data are defined as before. Analogously, the correspondence $L_\Theta$ preserves the tempered property, and the discrete property as well. More precisely, if the tempered $\phi = \phi_{GL} \oplus \phi_0 \oplus \phi_{GL}^\vee$ factors through a Levi subgroup $\hat{M} = \prod_{i \in I} GL(n_i, \mathbb{C}) \times Sp(2n_0, \mathbb{C})$ of $\hat{H}$, where $\phi_{GL}$ and $\phi_0$ are the discrete parameters of $\prod_{i \in I} GL(n_i)$ and $\tilde{G}_{2n_0} = \widetilde{Sp}_{2n_0}$ respectively. Denote $\tilde{M} = \prod_{i \in I} GL(n_i) \times \tilde{G}_{2n_0}$, and choose an arbitrary parabolic subgroup $\tilde{P} = \tilde{M}N$ of $\tilde{G}$, then*

$$\Pi_\phi^{\tilde{G}} = \bigsqcup_{\tilde{\sigma} \in \Pi_\phi^{\tilde{M}}} \{\text{the irreducible constituents of } I_{\tilde{P}}(\tilde{\sigma})\}.$$

*Furthermore, instead of discussing the properties of R-groups for $\tilde{G}$, we have the multiplicity-one property of the normalized induced representations from discrete series, i.e. for $\tilde{\sigma}_1, \tilde{\sigma}_2 \in \Pi_\phi^{\tilde{M}}$,*

- *$I_{\tilde{P}}(\tilde{\sigma}_i)$ is multiplicity-free for $i = 1, 2$;*
- *$\#J.H.(I_{\tilde{P}}(\tilde{\sigma}_1)) = \#J.H.(I_{\tilde{P}}(\tilde{\sigma}_2))$, here J.H. stands for the set of Jordan-Hölder constituents.*

*Proof.* One may consult [GS12, Theorem 1.3] for details of the theta correspondence for $(\tilde{G}, SO_{2n+1})$, and [Gol94, Art13, CG16] for the multiplicity-one property and cardinality relation for $SO_{2n+1}$. As for the archimedean case, the multiplicity-one statement follows from [Liar, Remark 7.4.4], and the constituents cardinality equality results from [AB98, Proposition 12.24]. □

2.2. **LLC for $\tilde{G}(\mathbb{R})$: Fine version.** In this subsection, we first recall the explicit local Langlands correspondence for $SO_{2n+1}(\mathbb{R})$ à la J. Adams (see [Ada10, Ada98] for details), and then verify the compatibility with the LLC for discrete series à la Shelstad and Arthur. As in [Ada10], we start with a complex group $\mathbb{H} = SO_{2n+1}(\mathbb{C})$ which is the isometry group of the quadratic form

$$q_0 = (-1)^n \begin{pmatrix} 1_{2n} & \\ & 1 \end{pmatrix}$$

of discriminant 1. Fix a Cartan subgroup $\mathbb{T} \subset \mathbb{H}$ with

$$\mathbb{T} = \underbrace{SO_2(\mathbb{C}) \times \cdots \times SO_2(\mathbb{C})}_{n \text{ copies}} \times 1$$

so that

$$X^*(\mathbb{T}) = \mathbb{Z}e_1 + \cdots + \mathbb{Z}e_n,\ X_*(\mathbb{T}) = \mathbb{Z}e_1^\vee + \cdots + \mathbb{Z}e_n^\vee,$$

where $e_i : \mathbb{T} \longrightarrow \mathbb{C}^\times$ is non-trivial only on the $i$-th copy of $SO_2$, on which it is given by the isomorphism

$$\begin{pmatrix} x & y \\ \bar{y} & x \end{pmatrix} \longrightarrow x + iy.$$

Moreover, $\{e_i^\vee\}$ is the basis of $X_*(\mathbb{T})$ dual to $\{e_i\}$. By fixing a Borel subgroup $\mathbb{B} \supset \mathbb{T}$, we may assume the systems of the positive simple roots and coroots are given by:

$$\Delta = \{e_1 - e_2, e_2 - e_3, \cdots, e_{n-1} - e_n, e_n\}, \qquad \Delta^\vee = \{e_1^\vee - e_2^\vee, e_2^\vee - e_3^\vee, \cdots, e_{n-1}^\vee - e_n^\vee, 2e_n^\vee\}.$$

The Langlands dual group of $\mathbb{H}$ is $\hat{\mathbb{H}} = Sp(2n)(\mathbb{C})$. It is equipped with a pair $\hat{\mathbb{T}} \subset \hat{\mathbb{B}}$ with

$$X^*(\hat{T}) = X_*(\mathbb{T}), \quad and \quad X_*(\hat{T}) = X^*(\mathbb{T})$$

and a system of simple roots $\hat{\Delta} = \Delta^\vee$.

The complex conjugation action on matrices defines an action of $Gal(\mathbb{C}/\mathbb{R})$ on $\mathbb{H}$, which gives rise to the unique compact real form $H_c$ with maximal torus $T_c = \mathbb{T} \cap H_c$. The complex and real Weyl groups are equal and given by

$$W(\mathbb{H}, \mathbb{T}) = W(H_c, T_c) \simeq S_n \ltimes W_0$$

with $W_0 = (\mathbb{Z}/2\mathbb{Z})^n$. Here, $S_n$ acts on $X^*(\mathbb{T})$ by permuting the $e_i$'s and the generator of the $i$-th copy of $\mathbb{Z}/2\mathbb{Z}$ in $W_0$ sends $e_i$ to $-e_i$ and fixes $e_j$ for $j \neq i$.

Let us set
$$X := \mathbb{T}[2] = T_c[2] = \{x \in \mathbb{T} : x^2 = 1\}.$$
Then one has a canonical isomorphism

(1)
$$X \xleftarrow{\sim} X_*(\mathbb{T})/2X_*(\mathbb{T}) = \bigoplus_{i=1}^n \mathbb{Z}/2\mathbb{Z} e_i$$
$$\lambda(-1) \longleftrightarrow \lambda.$$

The Weyl group $W(\mathbb{H}, \mathbb{T})$ preserves $X$ and acts through its quotient $S_n$. The isomorphism in (1) is $W(\mathbb{H}, \mathbb{T})$-equivariant.

Now let us consider the (pure) real forms of $\mathbb{H}$. Since $Out(\mathbb{H}) = 1$, the set of (pure) real forms of $\mathbb{H}$ up to isomorphism is classified by $H^1(\mathbb{R}, \mathbb{H})$, which is known to be equal to
$$X_1 := X/W(\mathbb{H}, \mathbb{T}).$$
More precisely, each $x \in X$ gives a 1-cocycle
$$c_x(\sigma) = x \quad (\text{with } \sigma \in Gal(\mathbb{C}/\mathbb{R}) \text{ non-trivial})$$
and $c_x$ is cohomologous to $c_{x'}$ if and only if $x$ and $x'$ are in the same $W(\mathbb{H}, \mathbb{T})$-orbit. Hence, $|X_1| = n + 1$, with each orbit $[x]$ determined by the number of nonzero coefficients of $x$ as a $\mathbb{Z}/2\mathbb{Z}$-linear combination of the $e_i$'s, or equivalently, the number of $+1$ when $x$ is regarded as an element of
$$\mathbb{T}[2] = \prod_{i=1}^n \{\pm 1\} \subset \prod_{i=1}^n SO_2(\mathbb{C}).$$
For each $x \in X$, one has a modified complex conjugation action on $\mathbb{H}$ given by:
$$h \longmapsto x \bar{h} x^{-1},$$
that gives an associated real form
$$H_x = \mathbb{H}^{Ad(x) \circ (-)}$$
which is the isometry group of the quadratic form determined by
$$q_x = (-1)^n x$$
of discriminant 1. A maximal compact subgroup of $H_x$ by complexification given by:
$$K_x = \mathbb{H}^{Ad(x)}.$$
Note that $T_c \subset H_x$ since the adjoint action of $x$ on $\mathbb{T}$ is trivial. Moreover, the real Weyl group of $H_x$ is
$$W(K_x, \mathbb{T}) = W(H_x, T_c) = W(\mathbb{H}, \mathbb{T})^{Ad(x)} = Stab_{W(\mathbb{H}, \mathbb{T})}(x).$$
For example, if $x \in \mathbb{T}[2]$ has $p$ $+1$'s and $q$ $-1$'s as its entries (with $p + q = n$), then
$$H_x = \begin{cases} SO(2p+1, 2q), & \text{if } n \text{ is even,} \\ SO(2q, 2p+1), & \text{if } n \text{ is odd,} \end{cases}$$
and its real Weyl group is isomorphic to:
$$(S_p \ltimes (\mathbb{Z}/2\mathbb{Z})^p) \times (S_q \ltimes (\mathbb{Z}/2\mathbb{Z})^q).$$

Let $\mathfrak{h} = Lie(\mathbb{H})$ be the complex Lie algebra of $\mathbb{H}$.

**Definition 2.2.1.** *A representation of a real form of $\mathbb{H}$ is a pair $(x, \pi)$ where $x \in X$ and $\pi$ an $(\mathfrak{h}, K_x)$-module. Say that $(x, \pi)$ is equivalent to $(x', \pi')$ if there exists $h \in \mathbb{H}$ (actually $h \in N_\mathbb{H}(\mathbb{T})$) such that $x' = hxh^{-1}$ and*
$$\pi' \simeq \pi^h := \pi \circ Ad(h^{-1}).$$
*Write $[x, \pi]$ for the equivalence class of $(x, \pi)$.*

Now suppose we are given a discrete series L-parameter

$$\phi : W_{\mathbb{R}} \longrightarrow \hat{\mathbb{H}} = Sp_{2n}(\mathbb{C}).$$

We may write:

$$\phi = \phi_1 \oplus \cdots \oplus \phi_n$$

with $\phi_i = Ind_{W_{\mathbb{C}}}^{W_{\mathbb{R}}}((\frac{z}{\bar{z}})^{a_i})$, $a_i \in \mathbb{Z} + 1/2$, $a_i > 0$. We may further assume $a_1 > a_2 > \cdots > a_n > 0$. The centralizer group $S_\phi = Z_{\hat{\mathbb{H}}}(\phi)$ is equal to:

$$\begin{aligned} S_\phi = \hat{T}[2] &\xrightarrow{\sim} X_*(\hat{T})/2X_*(\hat{T}) = X^*(\mathbb{T})/2X^*(\mathbb{T}) \\ \lambda(-1) &\longleftrightarrow \lambda \end{aligned} \quad (2)$$

Hence $S_\phi$ and $X$ are in canonical duality in view of (1) and (2). In particular, we have a canonical isomorphism

$$X \xrightarrow{\sim} \hat{S}_\phi = Irr(S_\phi). \quad (3)$$

For each $x \in X$, the work of Harish-Chandra furnishes an L-packet of discrete series representations $\Pi_x(\phi)$ of $H_x$. More precisely, set

$$\lambda_0 = \sum_{i=1}^n a_i x_i \in X^*(\mathbb{T}) \bigotimes_{\mathbb{Z}} \mathbb{C} = (Lie\ \mathbb{T})^*.$$

Then $\Pi_x(\phi)$ consists of all discrete series representations of $H_x$ with infinitesimal character $\lambda_0$. Indeed, Harish-Chandra associated to each $\lambda \in W(\mathbb{H}, \mathbb{T}).\lambda_0$ an irreducible discrete series representation $\pi_x(\lambda)$, such that $\pi_x(\lambda) \simeq \pi_x(\lambda')$ if and only if $\lambda \in W(K_x, \mathbb{T}).\lambda_0$. Hence

$$\Pi_x(\phi) = \{\pi_x(w^{-1}\lambda_0) : w \in W(\mathbb{H}, \mathbb{T})/W(K_x, \mathbb{T})\},$$

and the $W(K_x, \mathbb{T})$-orbit of $\lambda$ is the Harish-Chandra parameter of $\pi_x(\lambda)$.

Following D. Vogan, let us define the Vogan L-packet

$$\Pi(\phi) = \{[x, \pi_x(\lambda)] : x \in X, \lambda \in W(\mathbb{H}, \mathbb{T}).\lambda_0\}$$

where we have disregarded repetitions. If we pick a representative $x \in X$ for each $W(\mathbb{H}, \mathbb{T})$-orbit $\bar{x} \in X_1$, then we have bijections

$$\Pi(\phi) = \bigsqcup_{\bar{x} \in X_1} \Pi_x(\phi) \longleftrightarrow \bigsqcup_{\bar{x} \in X_1} W(\mathbb{H}, \mathbb{T})/W(K_x, \mathbb{T}). \quad (4)$$

We shall call this the <u>Harish-Chandra parametrization</u>, noting that it depends on the choice of representatives $x \in X$ for $\bar{x} \in X_1$. If the orbit $\bar{x}$ consists of elements of $X = \mathbb{T}[2] = \mu_2^n$ with $p$ +1-entries and $q$ -1-entries (with p+q=n), then we shall write the Harish-Chandra parameters of $\pi_x(\lambda) \in \Pi_x(\phi)$ as

$$\lambda = (a_{i_1}, \cdots, a_{i_p}; a_{j_1}, \cdots, a_{j_q})$$

where $\{a_{i_1}, \cdots, a_{i_p}, a_{j_1}, \cdots, a_{j_q}\}$ is a permutation of $\{a_1, \cdots, a_n\}$, with $a_{i_1} > \cdots > a_{i_p}, a_{j_1} > \cdots > a_{j_q}$.

J. Adams observed that one can give another parametrization of $\Pi(\phi)$ via the bijection

$$\begin{aligned} \Pi(\phi) &\longleftrightarrow X = \mathbb{T}[2] = \mu_2^n \\ [x, \pi_x(\lambda_0)] &\longleftrightarrow x \end{aligned} \quad (5)$$

We shall call this <u>Adams parametrization</u>. In view of the isomorphism (3), this gives a bijection

$$\Pi(\phi) \xrightarrow{naive} \hat{S}_\phi. \quad (6)$$

Under this bijection, the trivial character of $S_\phi$ corresponds to the irreducible representation of the compact group $H_c$ with highest weight $\lambda_0$. On the other hand, if

$$x_b = ((-1)^n, (-1)^{n-1}, \cdots, +1, -1) \in \mathbb{T}[2] = \mu_2^n,$$

then under (5), $x_b$ corresponds to the unique generic representation in $\Pi(\phi)$. Since the desiderate of the LLC requires the trivial character of $S_\phi$ to correspond to a generic representation, this suggests that one should modify the naive bijection (5) (and hence (6)) by translation-by $x_b$. Namely, we have a modified bijection

$$\Pi(\phi) \longleftrightarrow X \longleftrightarrow \hat{S}_\phi$$
(7)
$$[x, \pi_x(\lambda_0)] \longleftrightarrow x \cdot x_b \longleftrightarrow \chi_{xx_b}.$$

This is <u>Adams' version of the LLC</u>; we shall see shortly that it agrees with the LLC of Shelstad and Arthur.

The LLC of Shelstad and Arthur uses the Harish-Chandra parametrization (4) as its starting point. Recall that the bijection (4) requires the choice of a representative of $W(\mathbb{H}, \mathbb{T})$-orbits in $X$. For the orbit associated to the split group, we use the element $x_b \in X$, and the generic representation $\pi_{x_b}(\lambda)$ as the base point in $\Pi_{x_b}(\phi)$, so that

(8) $$\Pi_{x_b}(\phi) \longleftrightarrow W(\mathbb{H}, \mathbb{T})/W(K_x, \mathbb{T}).$$

We have an injection

(9)
$$W(\mathbb{H}, \mathbb{T})/W(K_{x_b}, \mathbb{T}) \hookrightarrow H^1(\mathbb{R}, \mathbb{T})$$
$$w \longleftrightarrow [c_w : \sigma \longmapsto \sigma_{x_b}(w).w^{-1}]$$

Note that:

$$H^1(\mathbb{R}, \mathbb{T}) = H^1(\mathbb{R}, X_*(\mathbb{T}) \bigotimes_{\mathbb{Z}} \mathbb{C}^\times)$$
(10)
$$= \prod_{i=1}^n H^1(\mathbb{R}, e_i^\vee \bigotimes_{\mathbb{Z}} \mathbb{C}^\times)$$
$$= \prod_{i=1}^n e_i \otimes_{\mathbb{Z}} \pi_0(\mathbb{R}^\times)$$
$$= \mathbb{T}[2] = X,$$

with $x \in X$ giving the 1-cocycle $\sigma \longmapsto x$. Thus, (8), (9) and (10) give

(11) $$\Pi_{x_b}(\phi) \hookrightarrow X.$$

For any other orbit $\bar{x} \in X_1$ with representative $x \in X$, the complex conjugation action on $\mathbb{H}$ associated to $H_x$ differs from that for $H_{x_b}$ by conjugation-by $xx_b$. This gives a 1-cocycle $c_x : \sigma \longmapsto x \cdot x_b$ with the associated class $[c_x] \in H^1(\mathbb{R}, \mathbb{T})$. Then we have:

(12)
$$\Pi_x(\phi) \longleftrightarrow W(\mathbb{H}, \mathbb{T})/W(K_x, \mathbb{T}) \hookrightarrow H^1(\mathbb{R}, \mathbb{T})$$
$$w \longleftrightarrow [\sigma \longmapsto \sigma_{x_b}(w)c_x(\sigma)w^{-1}]$$

Hence, we have a map (bijective)

(13) $$\Pi(\phi) = \bigsqcup_{\bar{x} \in X_1} \Pi_x(\phi) \longrightarrow H^1(\mathbb{R}, \mathbb{T}) = X.$$

This map is in fact easily described by:

**Lemma 2.2.2.** *The map given in (13) is given by:*

$$[x, \pi_x(\lambda_0)] \longmapsto [\sigma \longmapsto xx_b]$$

*In particular, it is independent of the choice of the orbit representatives $x \in X$ for $\bar{x} \in X_1$.*

*Proof.* Let us consider the embedding (11) for the split group. If $x = w.x_b$ for $w \in S_n \leq W(\mathbb{H}, \mathbb{T})$, then the image of $[x, \pi_x(\lambda_0)] = [x_b, \pi_{x_b}(w^{-1}\lambda_0)]$ under (11) is by definition the 1-cocycle $\sigma \longmapsto \sigma_{x_b}(w)w^{-1}$. But

$$\sigma_{x_b}(w)w^{-1} = x_b(wx_bw^{-1}) = x_bx.$$

So the 1-cocycle in question is $\sigma \longmapsto xx_b$. The case for the other orbits is similarly treated. Indeed, for the representative $x \in X$, its image under (12) is by definition the 1-cocycle $c_x : \sigma \longmapsto xx_b$. For any $w \in W(\mathbb{H}, \mathbb{T})$, the image of $[x, \pi_x(w^{-1}\lambda_0)] = [wx, \pi_{wx}(\lambda_0)]$ is by definition the 1-cocycle

$$\sigma \longmapsto \sigma_{x_b}(w)c_x(\sigma).w^{-1} = x_bwx_b.x.x_b.w^{-1} = x_b.wx.$$

Thus proves the Lemma. $\square$

As a consequence of the Lemma, the map (13) agrees with the map in (7), i.e. the LLC of Shelstad–Arthur agrees with that of Adams. We summarize the above discussion in the following table.

| | Harish-Chandra | Adams | Langlands–Shelstad |
|---|---|---|---|
| Parameter set for $\Pi^{\mathbb{H}}(\phi)$ | $\bigsqcup_{x \in X/W_{\mathbb{C}}^{\mathbb{H}}} W_{\mathbb{C}}^{\mathbb{H}}/W_{\mathbb{C}}^{K_x}.\lambda_0$ | $X = \mu_2^n$ | $\hat{S}_\phi = Irr(\mu_2^n)$ |
| $\lambda$ | $(a_{i_1}, \cdots, a_{i_p}; a_{j_1}, \cdots, a_{j_q})$ | $x_{i_k} = 1, x_{j_k} = -1$ | $\chi_{i_k} = (-1)^{n+1-i_k}, \chi_{j_k} = (-1)^{n+2-j_k}$ |
| generic | $(\cdots, a_{n-3}, a_{n-1}; \cdots, a_{n-2}, a_n)$ | $x_i = (-1)^{n+1-i}$ | 1 |
| compact $\lambda_0$ | $(a_1, \cdots, a_n;)$ | 1 | $\chi_i = (-1)^{n+1-i}$ |

TABLE 1. LLC for $SO(2n+1)$

<u>LLC for $\tilde{G}$</u>. For $\tilde{G} = Mp_{2n}(\mathbb{R})$, the discrete series $L$-packet $\phi$ also gives rise to an $L$-packet $\tilde{\Pi}(\phi)$, consisting of those discrete series representations with infinitesimal character $\lambda_0$. We would like to parametrize the representations in $\tilde{\Pi}(\phi)$ using $\hat{S}_\phi \simeq X$. One way is to proceed as in the Harish-Chandra parametrization. If we fix a Whittaker datum for $\tilde{G}$, there is a unique representation $\tilde{\pi}_g \in \tilde{\Pi}(\phi)$ which is generic with respect to this Whittaker datum, let's say $\tilde{\pi}_g$ has Harish-Chandra parameter $\tilde{\lambda}_g$. Then one has bijections

$$\tilde{\pi}(\phi) \longleftrightarrow W(U(n), \mathbb{T}) \backslash W(Sp(2n), \mathbb{T}).\tilde{\lambda}_g \longleftrightarrow H^1(\mathbb{R}, \mathbb{T}) \simeq X \simeq \mu_2^n \quad (14)$$
$$w \longleftrightarrow w^{-1}\sigma(w).$$

Note that $W(Sp(2n), \mathbb{T}) = W(U(n), \mathbb{T}) \ltimes W_0 = S_n \ltimes (\mathbb{Z}/2\mathbb{Z})^n$. So we have a canonical set of coset representatives for $W(U(n), \mathbb{T}) \backslash W(Sp(2n), \mathbb{T})$, namely $W_0 = (\mathbb{Z}/2\mathbb{Z})^n$.

On the other hand, one could use the theta correspondence to transport Shelstad's LLC for $SO_{2n+1}$ to $Mp_{2n}$. Let us fix an additive character $\psi_a(x) = e^{2\pi i a x}$ of $\mathbb{R}$ with $a > 0$. In [Ada98], J. Adams showed that the theta correspondence for $(SO(p,q), Mp(2n)(\mathbb{R}))$ with respect to $\psi_a$ gives a bijection

$$X \xrightarrow{LLC} \Pi(\phi) \xrightarrow{\Theta_{\psi_a}} \tilde{\Pi}(\phi) \quad (15)$$
$$x \longleftrightarrow [xx_b, \pi_{xx_b}(\lambda_0)]$$

For $x \in X$ such that $[x, \pi_x(\lambda_0)]$ has Harish-Chandra parameter $(a_{i_1}, \cdots, a_{i_p}; a_{j_1}, \cdots, a_{j_q})$ (so $x_{i_k} = 1$ and $x_{j_k} = -1$), $\Theta_{\psi_a}([x, \pi_x(\lambda_0)])$ has Harish-Chandra parameter

$$(a_{i_1}, \cdots, a_{i_p}; -a_{j_q}, \cdots, -a_{j_1}).$$

Moreover, $\Theta_{\psi_a}([x_b, \pi_{x_b}(\lambda_0)])$ is the unique $\psi_a$-generic representation in $\tilde{\Pi})(\phi)$. Now we have:

**Lemma 2.2.3.** *The two bijections $X \longleftrightarrow \tilde{\Pi}(\phi)$ in (14) and (15) are equal. (where we use the $\psi_a$-Whittaker datum in (14))*

*Proof.* Starting with $x \in X$, the bijection in (15) produces the representation $\Theta_{\psi_a}([xx_b, \pi_{xx_b}(\lambda_0)])$ which has Harish-Chandra parameter given by:

$$(a_i, i \in I^+_{xx_b}; -a_j, j \in I^-_{xx_b})$$

where

$$\begin{cases} I^+_{xx_b} = \{1 \le i \le n : (xx_b)_i = +1\} \\ I^-_{xx_b} = \{1 \le j \le n : (xx_b)_j = -1\} \end{cases}$$

In the bijection (14), the map

$$(\mathbb{Z}/2\mathbb{Z})^n \simeq W_0 \longrightarrow H^1(\mathbb{R}, \mathbb{T}) \simeq X = \mu_2^n$$

is the natural isomorphism (componentwise). If $x \in X$ gives rise to $w_x \in W_0$, then

$$w_x.\tilde{\lambda}_g = (a_i, i \in I^+_{xx_b}; -a_j, j \in I^-_{xx_b})$$

Since $\tilde{\lambda}_g$ is the Harish-Chandra parameter associated to $x_b$, i.e.
$$\tilde{\lambda}_g = (a_i, i \in I^+_{x_b}; a_j, j \in I^-_{x_b}).$$
□

We summarize this discussion in the following table.

| | Harish-Chandra | Adams | Langlands–Shelstad |
|---|---|---|---|
| Parameter set for $\Pi^{Mp_{2n}}_\phi$ | $W^G_\mathbb{C}/W^K_\mathbb{C} \simeq <\sigma_i>_{i=1,\cdots,n}$ | $X = \mu^n_2$ | $\hat{S}_\phi = Irr(\mu^n_2)$ |
| LLC | $\prod_{i_s \in I} \sigma_{i_s}, \quad I = \{i_1, \cdots, i_t\}$ | $x_{i_s} = -1, x_j = 1$ | $\chi_{i_s} = (-1)^{n+2-i_s}, \chi_j = (-1)^{n+1-j}$ |
| $\lambda_n$ | $(a_1, \cdots, a_n)$ | 1 | $\chi_i = (-1)^{n+1-i}$ |
| generic $\lambda_g$ | $(\cdots, a_{n-3}, a_{n-1}, -a_n, -a_{n-2}, \cdots)$ | $x_i = (-1)^{n+1-i}$ | 1 |
| general $\lambda$ | $(a_{i_1}, \cdots, a_{i_p}, -a_{j_q}, \cdots, -a_{j_1})$ | $x_{i_k} = 1, x_{j_k} = -1$ | $\chi_{i_k} = (-1)^{n+1-i_k}, \chi_{j_k} = (-1)^{n+2-j_k}$ |

TABLE 2. LLC for $\tilde{G}$

Here $\sigma_i$ is the reflection with respect to $a_i$, i.e. $\sigma_i(\cdots, a_i, \cdots) = (\cdots, -a_i, \cdots)$. And $j$ varies over the subset $\{1, \cdots, n\} \setminus \{i_1, \cdots, i_t\}$.

2.3. **Gan–Ichino multiplicity formula.** Given the local Langlands correspondence for $\tilde{G}$, we would like to recall Gan–Ichino's global multiplicity formula for elliptic tempered A-parameters of $\tilde{G}$ in $L^2_{disc}(G(\dot{F}) \backslash \tilde{G}(\mathbb{A}))$ (cf. [GI17]). Suppose that $\dot{F}$ is a number field, with $\mathbb{A}$ its associated Adèle ring. Fix an additive character of $\mathbb{A}/\dot{F}$ as in [GL16]. We investigate the multiplicity of any representation of $Mp_{2n}(\mathbb{A})$ in the $\Psi$-near equivalence subspace $L^2_\Psi(Mp_{2n})$ if $\Psi$ is an elliptic tempered A-parameter. As is well known, such $\Psi$ is a multiplicity-free sum of irreducible symplectic cuspidal automorphic representations $\Psi_i$ of $GL_{n_i}(\mathbb{A})$ as follows:
$$\Psi = \bigoplus_i \Psi_i,$$
where $\Psi_i$ is of symplectic type for all $i$. We may then formally define its centralizer as a free $\mathbb{Z}/2\mathbb{Z}$-module
$$S_\Psi = \bigoplus_i \mathbb{Z}/2\mathbb{Z}\, a_i$$
with a basis $\{a_i\}_i$, where each $a_i$ corresponds to $\Psi_i$. At any place $\nu$ of $\dot{F}$, this gives rise to a local L-parameter $\Psi_\nu$ together with a canonical map $S_\Psi \longrightarrow S_{\Psi_\nu} \longrightarrow A_{\Psi_\nu} := \pi_0(S_{\Psi_\nu})$. Then we have a compact group $A_{\Psi,\mathbb{A}} = \prod_\nu A_{\Psi_\nu}$ equipped with the diagonal map $\Delta : S_\Psi \longrightarrow A_{\Psi,\mathbb{A}}$. Note that for any $\eta = \otimes_\nu \eta_\nu \in \hat{A}_{\Psi,\mathbb{A}}$, we may form an irreducible genuine representation
$$\tilde{\pi}_\eta = \bigotimes_\nu \tilde{\pi}_{\eta_\nu}$$
of $Mp_{2n}(\mathbb{A})$, where $\tilde{\pi}_{\eta_\nu}$ is the associated genuine representation of $Mp_{2n}(\dot{F}_\nu)$ in the sense of the local Langlands correspondence of $\tilde{G}$ (see Section 2.1). Last but not least, we define a quadratic character $\epsilon_\Psi$ of $S_\Psi$ by setting
$$\epsilon_\Psi(a_i) = \epsilon(1/2, \Psi_i),$$
where $\epsilon(1/2, \Psi_i) \in \{\pm 1\}$ is the root number of $\Psi_i$.

**Theorem 2.3.1** (Global multiplicity formula). *(see [GI17]) Let $\Psi$ be an elliptic tempered A-parameter for $Mp_{2n}$. Then we have*
$$L^2_\Psi(Mp_{2n}) \simeq \bigoplus_{\eta \in \hat{A}_{\Psi,\mathbb{A}}} m_\eta \tilde{\pi}_\eta,$$
*where*
$$m_\eta = \begin{cases} 1, & \text{if } \Delta^*(\eta) = \epsilon_\Psi; \\ 0, & \text{otherwise.} \end{cases}$$
*In particular, $L^2_\Psi(Mp_{2n})$ is multiplicity-free.*

*Proof.* One may consult [GI17] for details. □

## 3. ENDOSCOPY THEORY AND TRACE FORMULA

In this section, we recall some results of the endoscopy theory which will be used in the proof of the endoscopic character identities. For full details, please refer to W.W. Li's papers [Li11, Li15, Liar].

3.1. **Endoscopy.** Let $F$ be a non-archimedean field of characteristic zero. Fix a symplectic $F$-vector space $(W, \langle \cdot, \cdot \rangle)$ of dimension $2n$. Write $G = Sp(W)$ and its associated metaplectic group $\tilde{G}$. Before we continue further, we first recall the parametrization of semisimple conjugacy classes as follows which would be involved in the definition of transfer factor. For details, please refer to [Wal01] or [Li11].

3.1.1. *Regular semisimple conjugacy classes.*
- $Sp(W)$ with dim $W = 2n$: the regular semisimple conjugacy classes are parametrized by the following data $\mathcal{O}(K/K^\#, x, c)$:
  - $(K, \tau)$ $2n$-dimensional étale $F$-algebra with involution $\tau$. Denote $K^\#$ the $\tau$-fixed étale subalgebra of $K$.
  - $x \in K^\times$ such that $\tau(x) = x^{-1}$, and $K = F[x]$.
  - $c \in K^\times$ with $\tau(c) = -c$.
- Split $SO(V, q)$ with dim $V = 2m+1$: the strongly regular semisimple conjugacy classes are parametrized by the following data $\mathcal{O}(K/K^\#, x, c)$:
  - $(K, \tau)$ $2m$-dimensional étale $F$-algebra with involution $\tau$. Denote $K^\#$ the $\tau$-fixed étale subalgebra of $K$.
  - $x \in K^\times$ such that $\tau(x) = x^{-1}$, and $K = F[x]$.
  - $c \in K^\times$ such that $\tau(c) = c$.
- $\mathcal{O}(K_1/K_1^\#, x_1, c_1)$ and $\mathcal{O}(K_2/K_2^\#, x_2, c_2)$ are equivalent if and only if there exists an $F$-algebra isomorphism and involution $\sigma : (K_1, \tau_1) \xrightarrow{\sim} (K_2, \tau_2)$ such that $\sigma(x_1) = x_2$ and $\sigma(c_1)c_2^{-1} \in N_{K_2/K_2^\#}(K_2^\times)$.

3.1.2. *Transfer factor.* For the convenience of readers, we summarize the properties of transfer factor as follows. For details, one may consult [Li11].
- (Elliptic endoscopic groups) An elliptic endoscopic datum for $\tilde{G}$ is an ordered pair $(n', n'') \in \mathbb{Z}_{\geq 0}^2$ verifying $n' + n'' = n$. The corresponding endoscopic group is
$$H := H_{n', n''} = SO(2n' + 1) \times SO(2n'' + 1).$$
- (Stable conjugacy) Two regular semisimple elements $\tilde{\delta}_1, \tilde{\delta}_2 \in \tilde{G}$ are called stably conjugate whenever $\delta_1$ and $\delta_2$ are stably conjugate in $G$, and $\Theta_\psi(-\tilde{\delta}_1) = \Theta_\psi(-\tilde{\delta}_2)$. Here $\Theta_\psi$ is the character of Weil representation $\omega_\psi$ of $\tilde{G}$, and $-\tilde{\delta}_i = (-1) \cdot \tilde{\delta}_i$ with $(-1)$ the canonical lift in $\tilde{G}$ of $-1 \in G$ (see [Li11, Definition 2.9]).
- (Norm correspondence) Fix an endoscopic datum $(n', n'')$. We say that two regular semisimple elements $\tilde{\delta} \in \tilde{G}$ and $\gamma = (\gamma', \gamma'') \in H(F)$ are of norm correspondence, if $\delta \in G(F)$ and $\gamma = (\gamma', \gamma'') \in H(F)$ correspond in terms of parameters:
$$\gamma \in \mathcal{O}(K'/K'^\#, x', c') \times \mathcal{O}(K''/K''^\#, x'', c'') \longleftrightarrow \delta \in \mathcal{O}(K/K^\#, (x', -x''), c) \text{ with } K = K' \times K''.$$
In such case, we say $(\tilde{\delta}, \gamma)$ is a norm pair.
- (Transfer factor) For a norm pair $(\gamma, \tilde{\delta})$ with
$$(\gamma, \delta) \in (\mathcal{O}(K'/K'^\#, x', c') \times \mathcal{O}(K''/K''^\#, x'', c''), \mathcal{O}(K/K^\#, (x', -x''), c)),$$
where $K = K' \times K''$. We define the transfer factor as follows.
$$\Delta(\gamma, \tilde{\delta}) := \frac{\Theta'_\psi}{|\Theta'_\psi|}(-\tilde{\delta}') \cdot \frac{\Theta''_\psi}{|\Theta''_\psi|}(\tilde{\delta}'') \cdot sgn_{K''/K''^\#}(P_{x'}(-x'')(x'')^{-n'} det(\delta' + 1)),$$
where $\Theta'_\psi$ (*resp.* $\Theta''_\psi$) is the Harish-Chandra character of the Weil representation $\omega'_\psi$ (*resp.* $\omega''_\psi$) of $\tilde{G}' := Mp(W') \times_{\mu_2} \mu_8$ (*resp.* $\tilde{G}''$), and $P_{x'} \in F[T]$ is the characteristic polynomial of $x' \in K'^\times$. Conventionally, $\Delta(\gamma, \tilde{\delta}) := 0$ if $(\gamma, \tilde{\delta})$ is not a norm correspondence pair. Such defined transfer factor has the following known properties:

- (*Genuine*) $\Delta(\gamma, \epsilon\tilde{\delta}) = \epsilon\Delta(\gamma, \tilde{\delta})$ for $\epsilon \in Ker(p)$.
- (*Cocycle property*) If $\tilde{\delta}$ and $\tilde{\delta}_1$ are stably conjugate, $\Delta(\gamma, \tilde{\delta}_1) = \langle \kappa, inv(\delta, \delta_1) \rangle \Delta(\gamma, \tilde{\delta})$, where $inv(\delta, \delta_1)$ is the associated cocycle in $H^1(F, G_\delta)$, and the endoscopic character $\kappa$ is defined as the product map of the second component of $H^1(F, G_\delta) = H^1(F, T') \times H^1(F, T'') = \mu_2^{I_1'^*} \times \mu_2^{I_2''^*}$, where $I_1'^*$ and $I_2''^*$ are constants determined by $\gamma = (\gamma', \gamma'')$.
- (*Symmetric*) $\Delta_{n',n''}((\gamma', \gamma''), \tilde{\delta}) = \Delta_{n'',n'}((\gamma'', \gamma'), -\tilde{\delta})$.
- (*Normalization à la Waldspurger*) Let $K = G(\mathcal{O}_F) \subset G(F)$ and $K_H = H(\mathcal{O}_F) \subset H$. Then for norm correspondence pairs $(\gamma, \delta) \in K \times K_H$, $\Delta(\gamma, \delta) = 1$ provided $(\gamma, \delta)$ are of regular reduction.
- (*Product formula*) Suppose $(\gamma, \delta) \in H_{G-reg}(\dot{F}) \times G(\dot{F})$ is a norm correspondence pair, and $\delta = (\tilde{\delta}_\nu)_\nu$ in $\tilde{G}(\mathbb{A})$, then
  * $\Delta_\nu(\gamma, \tilde{\delta}_\nu) = 1$ for almost all place $\nu$;
  * $\prod_\nu \Delta(\gamma, \tilde{\delta}_\nu) = 1$.
- (*Parabolic descent*) If a norm pair $(\gamma, \delta)$ lies in a Levi subgroup $M_H \times M$ with
  $$M_H = \prod_{i \in I}(GL(n_i') \times GL(n_i'')) \times H^b, \ M = \prod_{i \in I} GL(n_i) \times G^b, \text{ and } n_i = n_i' + n_i'',$$
  then $(\gamma^b, \delta^b)$ is also a norm pair in $H^b \times G^b$. Denote by $\Delta_{H,\tilde{G}}$ and $\Delta_{H^b,\tilde{G}^b}$ the transfer factors associated to $(H, G)$ and $(H^b, G^b)$ respectively, then
  $$\Delta_{H,\tilde{G}}(\gamma, \tilde{\delta}) = \Delta_{H^b, \tilde{G}^b}(\gamma^b, \tilde{\delta}^b),$$
  where $\tilde{\delta}^b$ is given by the relation: $j(\sigma_{GL}(\delta_{GL}), \tilde{\delta}^b) = \tilde{\delta}$, with $\delta = \delta_{GL} \times \delta^b$ and $\sigma_{GL}$ the natural splitting defined in [Luoar, Section 2.1]. Here $j : \tilde{G}' \times \tilde{G}'' \longrightarrow \tilde{G}$ is the natural product map with respect to the decomposition $W = W' + W''$.

3.1.3. *Transfer of orbital integrals.* For $x \in G$, set $G_x := C_G(x)^o$, and let
$$D_G(x) = |det(1 - Ad(x)|_{Lie\ G/Lie\ G_x})|^{1/2}.$$

Let $C_c^\infty(\tilde{G})_{--}$ be the anti-genuine subspace of $C_c^\infty(\tilde{G})$, i.e. $\tilde{f}(\epsilon\tilde{x}) = \epsilon^{-1}\tilde{f}(\tilde{x})$, similarly for other groups and function spaces. We define the normalized (stable) orbital integral on $\gamma \in H_{reg}$ for $f \in C_c^\infty(H(F))$ as

(H) $\quad O_\gamma(f) = D_H(\gamma) \int_{H_\gamma(F) \backslash H(F)} f(h^{-1}\gamma h)dh, \quad SO_\gamma(f) = D_H(\gamma) \int_{(H_\gamma \backslash H)(F)} f(h^{-1}\gamma h)dh.$

Similarly, for $\tilde{\delta} \in \tilde{G}_{reg}$, and $\tilde{f} \in C_c^\infty(\tilde{G})_{--}$,

(G) $\quad O_{\tilde{\delta}}(\tilde{f}) = D_G(\delta) \int_{G_\delta(F) \backslash G(F)} \tilde{f}(\tilde{g}^{-1}\tilde{\delta}\tilde{g})dg$

*Remark* 1. *We use the unique Haar measures on $G$ and $H$ such that the measures of $G(\mathcal{O}_F)$ and $H(\mathcal{O}_F)$ are both 1, and the Haar measures for* (H) *and* (G) *are compatible i.e. they are defined via the canonical isomorphisms between the centralizers of regular elements.*

As in [Liar], we set
$$\mathcal{I}(\tilde{G})_{--} := \{O_{\tilde{\delta}}(\tilde{f}) : \tilde{f} \in C_c^\infty(\tilde{G})\},$$
$$S\mathcal{I}(H_{n',n''}) := \{SO_\gamma(f) : f \in C_c^\infty(H_{n',n''}(F))\};$$

$\mathcal{I}_{cusp}(\tilde{G})_{--} :=$ the subspace in $\mathcal{I}(\tilde{G})_{--}$ of elements supported on the elliptic set,

$S\mathcal{I}_{cusp}(H_{n',n''}(F)) :=$ the subspace in $S\mathcal{I}(H_{n',n''}(F))$ of elements supported on the elliptic set.

Note that these spaces consist of certain functions on semisimple regular conjugacy classes. The following big theorem is the culmination of the fundamental work of W.W. Li [Li11, Liar] plus an $\epsilon$-contribution of the author [Luoar].

**Theorem 3.1.1.**

(i) (Transfer Theorem) Given $\tilde{f} \in C_c^\infty(\tilde{G})_{--}$, there exists $f^H \in C_c^\infty(H(F))$ such that
$$\sum_{\tilde{\delta}} \Delta(\gamma, \tilde{\delta}) O_{\tilde{\delta}}(\tilde{f}) = SO_\gamma(f^H)$$
for all $\gamma \in H_{G-reg}(F)$. We say $(\tilde{f}, f^H)$ is a transfer pair for $(\tilde{G}, H(F))$.

(ii) (Fundamental Lemma) Suppose the residue characteristic $p$ is not a power of $2$, $K = G(\mathcal{O}_F) \subset G(F)$ and $K_H = H(\mathcal{O}_F) \subset H(F)$, we define $\mu_K(\epsilon x) := \epsilon^{-1}$ if $x \in K$, otherwise $0$. Then $(\mu_K, 1_{K_H})$ is a transfer pair provided $mes(K) = mes(K_H) = 1$.

More generally, there is a natural algebra homomorphism
$$b : \mathcal{H}(\tilde{G})_{--} \longrightarrow \mathcal{H}(H)$$
of (anti-genuine) spherical Hecke algebras such that for any $\tilde{f} \in \mathcal{H}(\tilde{G})_{--}$, one may take $f^H = b(\tilde{f}) \in \mathcal{H}(H)$.

(iii) (Isomorphism Theorem) The collective transfer map
$$\mathcal{T}^\mathcal{E} : \mathcal{I}_{cusp}(\tilde{G})_{--} \longrightarrow \bigoplus_{H_{n',n''}} S\mathcal{I}_{cusp}(H_{n',n''})$$
is an isomorphism.

Moreover, this isomorphism is an isometry under the natural elliptic inner product (see [Liar, Corollary 5.2.5]).

3.1.4. *Compatibility with parabolic descent.* Let $\tilde{P} = \tilde{M}N$ be a standard parabolic subgroup of $\tilde{G}$. Let $K = G(\mathcal{O}_F) \subset G(F)$. Denote $\tilde{K}$ to be the preimage of $p : \tilde{G} \longrightarrow G$. For $\tilde{f} \in C_c^\infty(G(F))_{--}$, define the constant term along $P$ as
$$\bar{\tilde{f}}^{(P)}(\tilde{m}) = \delta_P(m)^{1/2} \int_N \bar{\tilde{f}}(\tilde{m}n) dn,$$
where
$$\bar{\tilde{f}}(\tilde{g}) = \int_{\tilde{K}} \tilde{f}(\tilde{k}\tilde{g}\tilde{k}^{-1}) d\tilde{k}, \qquad \text{and} \qquad \delta_P(m) = |det(Ad(m)|_{Lie\ N})|.$$
Here we take the unique Haar measures on $N$ and $\tilde{K}$ such that $mes(N \cap K) = mes(\tilde{K}) = 1$.

Before ending this subsection, we recall the parabolic descent property of transfer pairing $(\tilde{f}, f^H)$ which would be used later on. Given a Levi subgroup $M^{s_0}$ of $H = SO(2n'+1) \times SO(2n''+1)$ as
$$M^{s_0} = \prod_{i \in I} GL(n_i) \times SO(2m'+1) \times SO(2m''+1),$$
where $s_0 = (1_{m'}, -1_{m''})$. There exists an associated Levi subgroup $M$ of $G$:
$$M = \prod_{i \in I} GL(n_i) \times Sp(2m' + 2m'').$$
In such case, we denote $(I', I'')$ to be a partition of $I$ such that $n' = m' + \sum_{i \in I'} n_i$ and $n'' = m'' + \sum_{j \in I''} n_j$, and $M^{s_0} = (\prod_{i \in I'} GL(n_i) \times SO(2m'+1)) \times (\prod_{i \in I''} GL(n_i) \times SO(2m''+1))$ as a natural Levi subgroup of $H$.

**Definition.** Let $z[s_0] := ((z_i)_{i \in I}, 1) \in M^{s_0}(F)$ be the element of order $2$ defined by
$$z_i = \begin{cases} 1, & if\ i \in I'; \\ -1, & if\ i \in I''. \end{cases}$$
For $f \in C_c^\infty(M^{s_0})$, we define $f_{s_0}$ to be the twisted function by $z[s_0]$ as follows.
$$f_{s_0}(m) := f(z[s_0]m).$$

**Lemma 3.1.2.** *([Liar, Theorem 3.4.6]) Let $\tilde{G}$, $\tilde{M}$, $M^{s_0}$, and $H$ be as in the previous definition. For every $\tilde{f} \in C_{c,--}^\infty(\tilde{G})$, let $f^H$ be a transfer of $\tilde{f}$ to $H$. Then $\bar{\tilde{f}}^{(\tilde{M})}$ and $(\overline{f^H})_{s_0}^{(M^{s_0})}$ is a transfer pair for the Levi subgroups $(\tilde{M}, M^{s_0})$.*

3.2. **Trace formulas.** In this subsection, we would like to recall Arthur's and W.W. Li's trace formulas which will play an essential role in our local-global arguments.

3.2.1. *Stable trace formula: elliptic regular terms.* Note that W.W. Li has established the stable trace formula for elliptic terms [Li15], but we only need to use the elliptic regular part. For simplicity, we herein just state this part. Let $\Gamma_{rel}(G(\dot{F}))$ be the set of representatives for the elliptic regular semisimple conjugacy classes in $G(\dot{F})$, $\Sigma_{G-rel}(H(\dot{F}))$ the set of representatives for the elliptic $G$-regular semisimple stable conjugacy classes in $H(\dot{F})$; similarly for other groups.

**Definition 3.2.1.**
- For $\tilde{f} \in C_c^\infty(\tilde{G}(\mathbb{A}))_{--}$, we define the elliptic regular part of the trace formula in [Li15] for $\tilde{f}$ by $T_{rel}^{\tilde{G}}(\tilde{f})$ as follows.

$$T_{rel}^{\tilde{G}} := \sum_{\delta \in \Gamma_{rel}(G(\tilde{f}))} \tau(G_\delta) O_\delta(\tilde{f})$$

where $O_\delta(\tilde{f}) = \prod_\nu O_{\tilde{\delta}_\nu}(\tilde{f}_\nu)$ with $\delta = (\tilde{\delta}_\nu)_\nu$ in $\tilde{G}(\mathbb{A})$, for $\tilde{f} = \otimes_\nu \tilde{f}_\nu$, and $\tau(G_\delta)$ is the Tamagawa number of $G_\delta$ which equals to 1.

- Suppose $H$ is an endoscopic groups of $\tilde{G}$. For $f^H \in C_c^\infty(H(\mathbb{A}))$, we define the stable analogue $ST_{G-rel}^H(f^H)$ for $H$ as follows.

$$ST_{G-rel}^H(f^H) := \tau(H) \sum_{\gamma \in \Sigma_{G-rel}(H(\tilde{f}))} SO_\gamma^H(f^H)$$

where $SO_\gamma^H(f^H) = \prod_\nu SO_\gamma^H(f_\nu^H)$ for $f^H = \otimes_\nu f_\nu^H$, and $\tau(H)$ is the Tamagawa number of $H$.

**Theorem 3.2.2** (Stable trace formula: elliptic regular terms)**.** *Suppose $\tilde{f} = \otimes_\nu \tilde{f}_\nu \in C_c^\infty(\tilde{G})_{--}$ and an adélic transfer function $f^H = \otimes_\nu f_\nu^H \in C_c^\infty(H(\mathbb{A}))$ is chosen for each given elliptic endoscopic group $H := H_{n',n''}$. Then we have*

$$T_{rel}^{\tilde{G}}(\tilde{f}) = \sum_{\substack{H := H_{n',n''} \\ n'+n''=n}} \iota(\tilde{G}, H) ST_{G-rel}^H(f^H),$$

*where $\iota(\tilde{G}, H) = 1/2$ if one of $n'$ and $n''$ is zero, and $\iota(\tilde{G}, H) = 1/4$ otherwise.*

*Proof.* One may refer to [Li15, Theorem 5.2.6] or [Luoar] for details. $\square$

3.2.2. *Simple trace formulas.* Our second input is the stable trace formulas of Arthur. Before stating Arthur's simple trace formula, we introduce some notions on test functions. Say $\tilde{f} = \otimes_v \tilde{f}_v \in C_c^\infty(\tilde{G}(\mathbb{A}))_{--}$ is supercuspidal at a finite place $v_0$ if $trace\ \tilde{\pi}_{v_0}(\tilde{f}_{v_0}) = 0$ for all irreducible genuine non-supercuspidal representations $\tilde{\pi}$; one has a similar notion for test functions in $C_c^\infty(H(\mathbb{A}))$. Note that supercuspidal functions are cuspidal in the sense of Section 3.1.3.

**Arthur's simple trace formula.** *([Art88, Corollary 7.3 & 7.4]) Consider test functions $f = \otimes_v f_v' \in C_c^\infty(H(\mathbb{A}))$ which satisfy the following conditions:*
- *At some finite place $v_0$, $f$ is supercuspidal.*
- *At another finite place $v_1$, $f_{v_1}$ is supported on the elliptic regular semisimple locus of $H_{v_1}$.*

*Then*

$$\sum_{\gamma \in H(\dot{F})_{rel}/conj} vol(H_\gamma(\dot{F}) \backslash H_\gamma(\mathbb{A})) \int_{H_\gamma(\mathbb{A}) \backslash H(\mathbb{A})} f(x^{-1}\gamma x) dx$$
$$= \sum_{\pi \subset L_{cusp}^2(H(\dot{F}) \backslash H(\mathbb{A}))} m(\pi) trace\ \pi(f)$$

*where $m(\pi)$ is the multiplicity of $\pi$ in $L_{cusp}^2(H(\dot{F}) \backslash H(\mathbb{A}))$.*

**Simple trace formula for $\tilde{G}$.** *([Li14b, Theorem 6.7]) Consider anti-genuine test functions $\tilde{f} = \otimes_v \tilde{f}_v \in C_c^\infty(\tilde{G}(\mathbb{A}))_{--}$ which satisfy the following conditions:*
- *At some finite place $v_0$, $\tilde{f}$ is supercuspidal.*
- *At another finite place $v_1$, $\tilde{f}_{v_1}$ is supported on the elliptic regular semisimple locus of $\tilde{G}_{v_1}$.*

*Then*

$$\sum_{\delta \in G(\dot{F})_{rel}/conj} vol(G_\delta(\dot{F})\backslash G_\delta(\mathbb{A})) \int_{G_\delta(\mathbb{A})\backslash G(\mathbb{A})} \tilde{f}(x^{-1}\delta x)dx$$

$$= \sum_{\tilde{\pi} \subset L^2_{cusp}(G(\dot{F})\backslash \tilde{G}(\mathbb{A}))} m(\pi) trace\, \tilde{\pi}(\tilde{f})$$

*where $m(\tilde{\pi})$ is the multiplicity of $\tilde{\pi}$ in $L^2_{cusp}(G(\dot{F})\backslash \tilde{G}(\mathbb{A}))$.*

Our last input is Arthur's stable multiplicity formula for split $SO(2n+1)$. But for our purpose, we only need the following simple forms. For full details, one may consult [Art13, Theorem 4.1.2]. In what follows, for endoscopic groups $L'$ of $L$, we use the same notation $S^{L'}_{disc,\Psi}$ as in [Art13, Section 3.3] to indicate the contribution of an A-parameter $\Psi$ of the reductive group $L$ in the stable distribution $S^{L'}_{disc}$.

**Theorem 3.2.3.** *(Arthur's stable multiplicity formula) Given an elliptic tempered A-parameter $\Psi$ of $L := SO(2n+1)$, we denote $S_\Psi$ to be the global component group associated to $\Psi$ in $\hat{L} = Sp(2n,\mathbb{C})$ which is a finite abelian 2-group. Then for any $f = \otimes_\nu f_\nu \in \mathcal{H}(SO(2n+1))$, we have*

$$S^L_{disc,\Psi}(f) = \frac{2}{|S_\Psi|} \prod_\nu \Theta_{\Psi_\nu}(f_\nu).$$

*Moreover, if $\Psi$ factors through an elliptic endoscopic group $L' := SO(2n'+1) \times SO(2n''+1)$ of $L$, we further have*

(i) *If $\Psi = \Psi_1 \oplus \Psi_2$ factors through $\hat{L}'$ in a unique way up to $\hat{L}'$-conjugation, then for test functions $f' = f'_1 \otimes f'_2$ with $f'_1 = \otimes_\nu f'_{1,\nu}$ and $f'_2 = \otimes_\nu f'_{2,\nu}$, we have:*

$$S^{L'}_{disc,\Psi}(f') = \frac{2^?}{|S_\Psi|} \prod_\nu \Theta_{\Psi_{1,\nu}}(f'_{1,\nu}) \times \prod_\nu \Theta_{\Psi_{2,\nu}}(f'_{2,\nu}) = S^{SO(2n'+1)}_{disc,\Psi_1}(f'_1) \times S^{SO(2n''+1)}_{disc,\Psi_2}(f'_2),$$

*where $? = 2$ if $n' \neq 0$ and $n'' \neq 0$, and $? = 1$ otherwise.*

(ii) *If $n' = n''$, and $\Psi = \Psi_1 \oplus \Psi_2$ factoring through $\hat{L}'$ in two ways up to $\hat{L}'$-conjugation, i.e. $\Psi = \Psi_1 \oplus \Psi_2$ and $\Psi^{op} = \Psi_2 \oplus \Psi_1$, then we have*

$$S^{L'}_{disc,\Psi}(f') = \frac{2^2}{|S_\Psi|} \prod_\nu \Theta_{\Psi_{1,\nu}}(f'_{1,\nu}) \times \prod_\nu \Theta_{\Psi_{2,\nu}}(f'_{2,\nu}) + \frac{2^2}{|S_\Psi|} \prod_\nu \Theta_{\Psi_{2,\nu}}(f'_{1,\nu}) \times \prod_\nu \Theta_{\Psi_{1,\nu}}(f'_{2,\nu}).$$

(iii) *In general, we denote by $\Psi(L',\Psi)$ the set of equivalence classes of A-parameter $\Phi$ of $L'$ which gives rise to the A-parameter $\Psi$ of $L$. Then we have*

$$S^{L'}_{disc,\Psi}(f') = \frac{2^?}{|S_\Psi|} \sum_{\Phi \in \Psi(L',\Psi)} \prod_\nu \Theta_{\Phi_\nu}(f'_\nu),$$

*where $? = 2$ if $n' \neq 0$ and $n'' \neq 0$, and $? = 1$ otherwise.*

Building upon the above trace formulas, we shall obtain a simple stable trace formula for $\tilde{G}$ as follows.

**Corollary 3.2.4.** *(Simple stable trace formula I) Fix the notaions as above. Take simple test functions $\tilde{f} = \otimes_\nu \tilde{f}_\nu \in C^\infty_c(\tilde{G}(\mathbb{A}))_{--}$ such that*

- *at a finite place $\nu_1$, $\tilde{f}_{\nu_1}$ is cuspidal,*
- *at another finite place $\nu_2$, $\tilde{f}_{\nu_2}$ is supercuspidal and*
- *at an extra finite place $\nu_3$, $\tilde{f}_{\nu_3}$ is supported on the regular semisimple, elliptic locus of $\tilde{G}_{\nu_3}$.*

*For such test functions, we obtain a spectral simple stable trace formula for $\tilde{G}$ with respect to an elliptic tempered A-parameter $\Psi$ of $\tilde{G}$ as follows.*

$$\sum_{\substack{\otimes_\nu \eta_\nu \in \tilde{\mathcal{A}}_{\Psi,\mathbb{A}}: \\ \Delta^*(\otimes_\nu \eta_\nu) = \epsilon_\Psi}} \prod_\nu \Theta_{\tilde{\pi}_{\eta_\nu}}(\tilde{f}_\nu) = \sum_{(n',n'')} \frac{1}{|S_\Psi|} \sum_{\Phi \in \Psi(H_{n',n''},\Psi)} \prod_\nu \Theta_{\Phi_\nu}(f^{H_{n',n''}}_\nu).$$

*Here $f^{H_{n',n''}}$ is the associated transfer function on $H_{n',n''}$.*

*Proof.* Given the simple test functions $\tilde{f}$, the simple trace formula of $\tilde{G}$ gives rise to

$$(\star') \qquad T_{rel}^{\tilde{G}}(\tilde{f}) = I_{disc}^{\tilde{G}}(\tilde{f}) := \sum_{\tilde{\pi} \subset L^2_{disc}(G(\dot{F}) \backslash \tilde{G}(\mathbb{A}))} m(\tilde{\pi}) trace \, \tilde{\pi}(\tilde{f}).$$

Similarly, the simple stable trace formula [Liar, Theorem 8.1.4] gives

$$ST_{G-rel}^{H_{n',n''}}(f^{H_{n',n''}}) = S_{disc}^{H_{n',n''}}(f^{H_{n',n''}}).$$

Thus the geometrical stable trace formula in Theorem 3.2.2 gives rise to a spectral stable trace formula (see [Liar, Theorem 8.1.4])

$$I_{disc}^{\tilde{G}}(\tilde{f}) = \sum_{(n',n'')} \iota(\tilde{G}, H_{n',n''}) S_{disc}^{H_{n',n''}}(f^{H_{n',n''}}).$$

On the other hand, following Arthur's standard procedure on extracting the $\Psi$-part of the stable trace formula (cf. [Art13, Chapter 3]), Gan–Ichino's global multiplicity formula for $\Psi$ (see Theorem 2.3.1) says that

$$(\star\star') \qquad I_{disc,\Psi}^{\tilde{G}}(\tilde{f}) = \sum_{\substack{\otimes_\nu \eta_\nu \in \check{A}_{\Psi,\mathbb{A}}: \\ \Delta^*(\otimes_\nu \eta_\nu) = \epsilon_\Psi}} \prod_\nu \Theta_{\tilde{\pi}_{\eta_\nu}}(\tilde{f}_\nu).$$

On the other hand, Arthur's stable multiplicity formula (see Theorem 3.2.3) says that,

$$S_{disc,\Psi}^{H_{n',n''}}(f^{H_{n',n''}}) = \frac{2^?}{|S_\Psi|} \sum_{\Phi \in \Psi(H_{n',n''}, \Psi)} \prod_\nu \Theta_{\Phi_\nu}(f_\nu^{H_{n',n''}}).$$

Here ? = 2 if $n'n'' \neq 0$, otherwise ? = 1. Whence the corollary holds. □

*Remark* **2.** *As pointed out by the referee, it is a priori possible for all supercuspidal test functions of $\tilde{G}$ to be transferred to zero for a particular $H$. To overcome this difficulty, we apply Gan–Ichino's result on the no-embedded eigenvalues property for $\Psi$ (see [GI17, Remark 1.2]). We may then weaken our supercuspidal condition on $\tilde{f}_{\nu_2}$ to be just cuspidal. This is the right version stated as follows which we would use later on.*

**Corollary 3.2.5.** *(Simple stable trace formula II) Fix the notations as above. Take simple test functions $\tilde{f} = \otimes_\nu \tilde{f}_\nu \in C_c^\infty(\tilde{G}(\mathbb{A}))_{--}$ such that*
- *at two finite places $\nu_1, \nu_2$, $\tilde{f}_{\nu_1}$ and $\tilde{f}_{\nu_2}$ are cuspidal, and*
- *at an extra finite place $\nu_3$, $\tilde{f}_{\nu_3}$ is supported on the regular semisimple, elliptic locus of $\tilde{G}_{\nu_3}$.*

*For such test functions, we obtain a spectral simple stable trace formula for $\tilde{G}$ with respect to an elliptic tempered A-parameter $\Psi$ of $\tilde{G}$ as follows.*

$$\sum_{\substack{\otimes_\nu \eta_\nu \in \check{A}_{\Psi,\mathbb{A}}: \\ \Delta^*(\otimes_\nu \eta_\nu) = \epsilon_\Psi}} \prod_\nu \Theta_{\tilde{\pi}_{\eta_\nu}}(\tilde{f}_\nu) = \sum_{(n',n'')} \frac{1}{|S_\Psi|} \sum_{\Phi \in \Psi(H_{n',n''}, \Psi)} \prod_\nu \Theta_{\Phi_\nu}(f_\nu^{H_{n',n''}}).$$

*Here $f^{H_{n',n''}}$ is the associated transfer function of $H_{n',n''}$.*

*Proof.* The only difference is that, under this weaker test function condition, Equation $(\star')$ is replaced by

$$(\star'') \qquad T_{ell}^{\tilde{G}}(\tilde{f}) = \sum_M |W(G,M)|^{-1} \sum_{s \in W(G,M)_{reg}} |det(s-1)_{\mathfrak{a}_M/\mathfrak{a}_G}|^{-1} trace \, (M_P(s,0) I_{\tilde{P},disc}(0, \tilde{f})),$$

where the first sum runs over associated classes of standard Levi subgroups of $G$. We refer the reader to [Li14b, Art05] for the precise definitions of other terms in the formula above. Applying Gan–Ichino's result on the no embedded eigenvalue property for $\Psi$ (see [GI17]), we know that only the term $M = G$ contributes to the $\Psi$-part, which in turn gives rise to the above Equation $(\star\star')$. Then the corollary follows from the same argument as in Corollary 3.2.4. □

## 4. Transfer of parabolically induced representations

Recall that $F$ is a local field of characteristic 0. In this section, we begin the proof of our Main Theorem by induction on the rank of $\tilde{G}_{2n} = \widetilde{Sp}^{(8)}(2n, F)$. Hence we fix the endoscopic group $H = SO(2n'+1) \times SO(2n''+1)$ of $\tilde{G}$, and assume that the Main Theorem has been shown for $\tilde{G}_{2m} = \widetilde{Sp}^{(8)}(2m, F)$ for arbitrary $m < n$. We suppose that $\phi_H : L_F \longrightarrow {}^L H$ is a tempered $L$-parameter of $H$ which is not a discrete series. We shall establish the Main Theorem for such $\phi_H$'s in this section. The case when $\phi_H$ is discrete series will be addressed in the next section.

### 4.1. Borel subgroup.
We first consider the case when
$$\phi_H : WD_F \longrightarrow \hat{T}_H(\mathbb{C}) \hookrightarrow \hat{H}(\mathbb{C})$$
factors through a maximal torus $\hat{T}_H(\mathbb{C})$ of $\hat{H}(\mathbb{C})$. In this case, the associated $L$-packet of $H$ consists of a principal series representation induced from a character of $T_H \coloneqq T' \times T''$ of $H = SO(2n'+1) \times SO(2n''+1)$. Given unitary characters $\alpha_1, \alpha_2, \cdots, \alpha_{n'}$ and $\beta_1, \beta_2, \cdots, \beta_{n''}$ of $F^\times$, let $\mu_H = \mu'_H \times \mu''_H \coloneqq (\alpha_1 \otimes \alpha_2 \otimes \cdots \otimes \alpha_{n'}) \times (\beta_1 \otimes \beta_2 \otimes \cdots \otimes \beta_{n''})$ denote the character
$$diag(a_1, \cdots, a_{n'}, 1, a_1^{-1}, \cdots, a_{n'}^{-1}) \times diag(b_1, \cdots, b_{n''}, 1, b_1^{-1}, \cdots, b_{n''}^{-1}) \longmapsto \prod_{i=1}^{n'} \alpha_i(a_i) \prod_{j=1}^{n''} \beta_j(b_j)$$
of the diagonal maximal torus $T_H$, and let
$$\pi \coloneqq I_{B_H}^H(\mu_H) \qquad \text{(normalized induction)}$$
be the associated principal series representation. Note that the associated character distribution is stable as normalized induction preserves stability. In particular, we may consider its endoscopic transfer $Tran_H^{\tilde{G}}(\Theta_\pi)$ to $\tilde{G} = \widetilde{Sp}^{(8)}(2n, F)$.

On the other hand, the norm correspondence of their split tori is
$$T'_H \times T''_H \longleftrightarrow T_G = T'_G \times T''_G$$
$$\gamma = (\gamma', 1, \gamma'^{-1}) \times (\gamma'', 1, \gamma''^{-1}) \longmapsto \delta = (\gamma', \gamma'^{-1}, -\gamma'', -\gamma''^{-1}).$$
We denote by $\mu_G = \mu'_G \times \mu''_G$ the corresponding character of $T_G$ determined by $\mu_H$ via:
$$\begin{cases} \mu'_G(\gamma', \gamma'^{-1}) = \mu'_H(\gamma', 1, \gamma'^{-1}) \\ \mu''_G(\gamma'', \gamma''^{-1}) = \mu''_H(\gamma'', 1, \gamma''^{-1}) \end{cases}$$
Hence, if $\gamma \longleftrightarrow \delta$ under the norm correspondence, we have
$$\mu_H(\gamma) = \mu''_G(-1) \mu_G(\delta).$$
Before going to the endoscopic transfer, we make an observation about the transfer factor of W.W. Li concerning the split torus (cf. [Li11]).

**Lemma 4.1.1.** *For each norm correspondence pair $(\gamma, \delta) \in H(F)_{reg} \times G(F)_{reg}$ with $(\gamma, \delta)$ lying in the split torus, we have:*
$$\Delta(\gamma, \tilde{\delta}) = \epsilon(\tilde{\delta}) \gamma(det(\delta), \psi)^{-1},$$
*where $\tilde{\delta} = (\delta, \epsilon(\tilde{\delta}))$ with $\epsilon(\tilde{\delta}) \in \mu_8$, and $\gamma(\cdot, \psi)$ is the relative Weil index (see [Kud96]).*

*Proof.* Note that we have a splitting section $T_G \longrightarrow \tilde{G}$ given by
$$\delta \longmapsto (\delta, \gamma(det(\delta), \psi)),$$
and via [Li11, Corollaire 4.7]
$$\Theta_\psi((\delta, \gamma(det(\delta), \psi))) \in \mathbb{R}_{>0},$$
where $\Theta_\psi$ is the Harish-Chandra character of the Weil representation associated to $\psi$. By [Li11, Proposition 5.11], we have
$$\Delta(\gamma, \tilde{\delta}) = \frac{\Theta_\psi}{|\Theta_\psi|}(\tilde{\delta}) \gamma_\psi(q[C_{\delta'}])^{-1}$$

as $\Delta_0(\gamma, \delta) = 1$. Here $\delta' = (\gamma', \gamma'^{-1})$ and $q[C_{\delta'}]$ is a quadratic form associated to the Cayley form $C_{\delta'}$ (see [Liar, P.544]). By [Li11, Theorem 4.16] and [Kud96, Lemma 4.1], we know that

$$\gamma_\psi(q[C_{\delta'}])^{-1} = \gamma_\psi((-1)^{n-1})\gamma_\psi((-1)^n) = \gamma_\psi(1)^2 \gamma(-1, \psi) = \gamma(-1, \psi)^{-1}\gamma(-1, \psi) = 1,$$

whence the lemma holds. $\square$

Consider the particular genuine representation

$$\tilde{\pi} := I_{\tilde{B}_G}^{\tilde{G}}(\chi_\psi \mu_G),$$

where $\chi_\psi$ is the genuine character of $\tilde{T}_G$ defined by

$$\chi_\psi(x, \epsilon) = \epsilon \gamma(det(x), \psi)^{-1}.$$

Now we have:

**Proposition 4.1.2.** *The Harish-Chandra characters $\Theta_\pi$ and $\Theta_{\tilde{\pi}}$ satisfy the following character relation for matching test functions $\tilde{f} \in C_c^\infty(\tilde{G})_{--}$ and $f^H \in C_c^\infty(H)$:*

$$\Theta_\pi(f^H) = \mu_G''(-1)\Theta_{\tilde{\pi}}(\tilde{f}) = \epsilon(1/2, \phi'', \psi)\Theta_{\tilde{\pi}}(\tilde{f}),$$

*where $\phi''$ is the L-parameter corresponding to $I_{T_H''}^{H''}(\mu_H'')$, and $\epsilon(1/2, \phi'', \psi)$ is independent of $\psi$. In particular, for unramified characters $\mu_H$,*

$$\Theta_\pi(f^H) = \Theta_{\tilde{\pi}}(\tilde{f}).$$

*Proof.* By the Weyl integration formula, we have, see [CG15, Section 5.8]) for the details,

$$\Theta_\pi(f^H) = \int_{T_H(F)} \mu_H(t) O_t(f^H) dt.$$

Note that, by Lemma 4.1.1, we have the transfer pairing identity

$$O_\gamma(f^H) = \Delta(\gamma, \tilde{\delta})O_{\tilde{\delta}}(\tilde{f}) = \epsilon(\tilde{\delta})\gamma(det(\delta), \psi)^{-1}O_{\tilde{\delta}}(\tilde{f}).$$

Therefore we have

$$\Theta_\pi(f^H) = \mu_G''(-1) \int_{T_G(F)} \chi_\psi \mu_G(\delta) O_{\tilde{\delta}}(\tilde{f}) d\delta,$$

whence the proposition follows from the computation of symplectic local root number (see [GGP12, Proposition 5.1]), i.e. $\epsilon(1/2, \phi'', \psi) = \mu_H''(-1) = \mu_G''(-1)$. $\square$

4.2. **Parabolic subgroups.** Now we consider the general case where

$$\phi_H : WD_F \longrightarrow \hat{M}_H(\mathbb{C}) \hookrightarrow \hat{H}(\mathbb{C}),$$

is a tempered L-parameter factoring through a Levi subgroup $\hat{M}_H(\mathbb{C})$ of $\hat{H}(\mathbb{C})$ such that $\phi_H \in \Phi_2(M_H)$. The Levi subgroup $M_H$ has the form

$$M_H \simeq (\prod_{i \in I'} GL(n_i) \times SO(2m'+1)) \times (\prod_{i \in I''} GL(n_i) \times SO(2m''+1)).$$

By composition with the embedding

$$\iota : \hat{H}(\mathbb{C}) \longrightarrow \hat{\tilde{G}} = Sp_{2n}(\mathbb{C}),$$

we have a tempered L-parameter $\phi = \iota \circ \phi_H$ of $\tilde{G}$ which factors through the Levi subgroup $\hat{M}$ of $\hat{\tilde{G}}$:

$$\phi : WD_F \xrightarrow{\phi_M} \hat{M}(\mathbb{C}) \longrightarrow \hat{\tilde{G}}(\mathbb{C})$$

where

$$M \simeq \prod_{i \in I} GL_{n_i}(F) \times \widetilde{Sp}_{2n_0}^{(8)}(F),$$

with

$$\begin{cases} I = I' \cup I'' \\ n_0 = m' + m''. \end{cases}$$

We may write $\phi = \phi_{GL} \times \phi_0$, where $\phi_{GL}$ denotes the GL-part of $\phi$ and $\phi_0$ denotes the $Sp_{2n_0}$-part. As a 2n-dimensional representation, we have:
$$\phi = \phi_{GL} \oplus \phi_0 \oplus \phi_{GL}^\vee.$$

Now observe that $M_H$ is an endoscopic group of $M$ which is determined by an element $s_0 \in S_{\phi_0} = Z_{Sp_{2n_0}}(\phi_0)$ with $s_0^2 = 1$ via:
$$\hat{M}_H = Z_{\hat{M}}(s_0).$$

Likewise, the endoscopic group $H$ of $\tilde{G}$ is associated to an element $s \in S_\phi$ with $s^2 = 1$. We shall sometimes write:
$$M_H = M^{s_0} \quad \text{and} \quad H = \tilde{G}^s.$$

The element $s$ is determined by the tuple $(s_0, I', I'')$ via:

(16) $$s = (s_0, \prod_{i \in I'} 1_{GL_{n_i}}, \prod_{i \in I''} -1_{GL_{n_i}}) \in \hat{M}(\mathbb{C}) = Sp_{2n_0}(\mathbb{C}) \times \prod_{i \in I'} GL_{n_i}(\mathbb{C}) \times \prod_{i \in I''} GL_{n_i}(\mathbb{C}).$$

Moreover, there is a natural commutative diagram:

(17)
$$\begin{array}{ccc} s_0 \in S_{\phi_M} = S_{\phi_0} & \xrightarrow{\iota} & S_\phi \ni s \\ \downarrow & & \downarrow \\ \bar{s}_0 \in \pi_0(S_{\phi_0}) = A_{\phi_0} & \xrightarrow{\bar{\iota}} & A_\phi \ni \bar{s} \end{array}$$

Observe that, in general, $s \neq \iota(s_0)$, but one has

(18) $$\bar{s} = \bar{\iota}(\bar{s}_0) \in A_\phi,$$

where $\bar{s}$ and $\bar{s}_0$ are the images of $s$ and $s_0$ in $A_\phi$ and $A_{\phi_0}$ respectively.

We now consider the L-packets associated to $\phi_H$ and $\phi$. By our discussion in Section 2.1, we have:
$$\Pi_{\phi_H}^H = \bigsqcup_{\sigma \in \Pi_\phi^{M^{s_0}}} \{\text{irreducible constituents of } I_{M^{s_0}}^H(\sigma)\},$$

and
$$\Pi_\phi^{\tilde{G}} = \bigsqcup_{\tilde{\sigma} \in \Pi_\phi^{\tilde{M}}} \{\text{irreducible constituents of } I_{\tilde{M}}^{\tilde{G}}(\tilde{\sigma})\}.$$

We have seen before (see Section 2.1) that the above induced representations are multiplicity-free. Hence, we have:

($\star$) $$\Theta_{\phi_H}^H := \sum_{\pi \in \Pi_{\phi_H}^H} \Theta_\pi = \sum_{\sigma \in \Pi_{\phi_H}^{M^{s_0}}} \Theta_{I_{M^{s_0}}^H(\sigma)}.$$

Likewise, in view of (17) and (18) above, we have:

($\star\star$)
$$\sum_{\eta \in \hat{A}_\phi} \eta(s) \Theta_{\tilde{\pi}_\eta}$$
$$= \sum_{\eta^{\tilde{M}} \in \hat{A}_{\phi_{\tilde{M}}}} \sum_{\substack{\eta \in \hat{A}_\phi \\ \eta \circ \iota = \eta^{\tilde{M}}}} \eta(\iota(s_0)) \Theta_{\tilde{\pi}_\eta} \text{ by (2)}$$
$$= \sum_{\eta^{\tilde{M}} \in \hat{A}_{\phi_{\tilde{M}}}} \eta^{\tilde{M}}(s_0) \left( \sum_{\substack{\eta \in \hat{A}_\phi \\ \eta \circ \iota = \eta^{\tilde{M}}}} \Theta_{\tilde{\pi}_\eta} \right)$$
$$= \sum_{\eta^{\tilde{M}} \in \hat{A}_{\phi_{\tilde{M}}}} \eta^{\tilde{M}}(s_0) \Theta_{I_{\tilde{M}}^{\tilde{G}}(\tilde{\sigma}_{\eta^{\tilde{M}}})}$$

where the last equality follows from the fact (see [GS12, Theorem 1.3 (ii)])
$$\{\text{irreducible consituents of } I_{\tilde{M}}^{\tilde{G}}(\tilde{\sigma}_{\eta^{\tilde{M}}})\} = \{\tilde{\pi}_\eta : \eta \circ \iota = \eta^{\tilde{M}}\}.$$

Since the character distribution $\Theta^H_{\phi_H}$ is stable, we may consider the transfer lift of $\Theta^H_{\phi_H}$ to $\tilde{G}$. In view of the parabolic descent property of transfer lifting, i.e. Lemma 3.1.2, we may apply induction argument to deduce our desired endoscopic character identities for tempered packets. We assume first the expected endoscopic character identities (i.e. Main Theorem) hold for any discrete series packet of any $M^{s_0}$ and the associated tempered packet for $\tilde{M}_0 = \widetilde{Sp}^{(8)}_{2n_0}(F)$.

Now we are ready to prove the desired endoscopic character identities.

**Proposition 4.2.1.** *The following character identity holds for matching local test functions $\tilde{f}$ and $f^H$:*
$$\Theta^H_{\phi_H}(f^H) = \epsilon(1/2, \phi'', \psi) \sum_{\eta \in \hat{A}_\phi} \eta(s) \Theta_{\tilde{\pi}_\eta}(\tilde{f}).$$

*In particular, we have the expected endoscopic character lifting identity:*
$$Trans_{n',n''}(\Theta^H_{\phi_H}) = \epsilon(1/2, \phi'', \psi) \sum_{\eta \in \hat{A}_\phi} \eta(s) \Theta_{\tilde{\pi}_\eta}.$$

*Proof.* By the Weyl integration formula and Formula $(\star)$, we have
$$\Theta^H_{\phi_H}(f^H) = \sum_{T_{M^{s_0}}} \frac{1}{|W(M^{s_0}(F), T_{M^{s_0}}(F))|} \int_{T_{M^{s_0}}} \Theta^{M^{s_0}}_{\phi_H}(t) O_t(f^H) dt,$$

where the sum is over representatives of the conjugacy classes of the maximal tori of $M^{s_0}$. By the parabolic descent formula of orbital integral, we have:
$$\Theta^H_{\phi_H}(f^H) = \Theta^{M^{s_0}}_{\phi_H}((\overline{f^H})^{(M^{s_0})}).$$

Based on Formula $(\star\star)$, the same argument implies
$$\sum_{\eta \in \hat{A}_\phi} \eta(s) \Theta_{\tilde{\pi}_\eta}(\tilde{f}) = \sum_{\eta \in \hat{A}_{\phi_{\tilde{M}}}} \eta(s_0) \Theta^{\tilde{M}}_{\tilde{\pi}_\eta}(\bar{\tilde{f}}^{(\tilde{M})}),$$

Note that, by Lemma 3.1.2, the GL-part transfer factor is twisted by $z[s]$, and readily the resulting effect is that the character identities of GL-part should be twisted by the associated central sign, thus by assumption we have:
$$\sum_{\eta \in \hat{A}_\phi} \eta(s) \Theta_{\tilde{\pi}_\eta}(\tilde{f}) = \epsilon(1/2, \phi''_{GL}, \psi)^{-1} \epsilon(1/2, \phi''_0, \psi)^{-1} \Theta^{M^{s_0}}_{\phi_H}((\overline{f^H})^{(M^{s_0})}) = \epsilon(1/2, \phi'', \psi)^{-1} \Theta^H_{\phi_H}(f^H),$$

where the last identity follows from the additive property of local root number (cf. [GP92, P.18]). Whence the proposition holds. □

## 5. Transfer of discrete series characters: archimedean case

In this section, we shall reinterpret D. Renard's endoscopic character identities for $\tilde{G}(\mathbb{R})$ to fit into the framework of the Main Theorem. Note that W.W. Li has reformulated D. Renard's work to fit into his endoscopy theory for $\tilde{G}$ (cf. [Liar]). Thus what we should do is to express W.W. Li's transfer identity in terms of the characters of component groups.

We first suppose $F = \mathbb{R}$. Fix a non-trivial character $\psi = \psi_a = e^{2\pi i a x} : F \longrightarrow S^1$ with $a > 0$, $G = Sp(2n)$ and $\tilde{G} \coloneqq \widetilde{Sp}^{(8)}_\psi(2n) \simeq Mp(2n) \times_{\mu_2} \mu_8$. Fix a maximal compact subgroup $K$ of $G(\mathbb{R})$ and an anisotropic maximal torus $T \subset G(\mathbb{R})$ such that $T \subset K$. Consider an endoscopic group $H \coloneqq SO(2p_0 + 1) \times SO(2q_0 + 1)$ of $\tilde{G}$ with $p_0 + q_0 = n$. For discrete series L-parameter of $SO(p,q)$ associated to a Harish-Chandra parameter
$$\lambda_p = (a_1, \cdots, a_{p_0}; b_1, \cdots, b_{q_0})$$

with $p_0 = [p/2]$ and $q_0 = [q/2]$, where $a_i, b_j \in \mathbb{Z} + 1/2$, $a_1 > \cdots > a_{p_0} > 0, b_1 > \cdots > b_{q_0} > 0$, then the theta-lift with respect to $\psi$ of $\pi_{\lambda_p}$ to $Mp(2n)(\mathbb{R})$ is the discrete series representation $\tilde{\pi}_{\lambda'_p}$ with Harish-Chandra parameter
$$\lambda'_p = (a_1, \cdots, a_{p_0}, -b_{q_0}, \cdots, -b_1).$$

Notice that we may also regard $\lambda_p$ as Harish-Chandra parameter of $H$, i.e. $(a_1, \cdots, a_{p_0}) \times (b_1, \cdots, b_{q_0})$. For this pair $(\lambda_p, \lambda'_p)$, we denote by $\phi$ the corresponding L-parameter, and $\Theta^H_\phi = \Theta^H_{\lambda_p}$ the stable distribution of $H$ associated to $\lambda_p$. Then we have the following theorem which was proved by D. Renard (cf. [Ren99, Proposition 6.2]) and later reformulated by W.W. Li (cf. [Liar, Theorem 7.3.7]).

**Theorem 5.0.1.** *Use the same notations defined above, we then have:*

$$Tran_{p_0,q_0}(\Theta^H_{\lambda_p}) = e^{\lambda'_p - \rho^\diamond}((1_{p_0}, -1_{q_0})) \sum_{w \in W^G_\mathbb{C}/W^K_\mathbb{C}} \kappa_T(w^{-1}) \Theta^{\tilde{G}}_{\tilde{\pi}_{w^{-1}\lambda'_p}},$$

where $\rho^\diamond = (p_0, \cdots, 1) \times (q_0, \cdots, 1)$, $e^{\lambda'_p - \rho^\diamond}$ *is a genuine character of $\tilde{T}$ (cf. [Liar, Ren99]) and $\kappa_T$ is the endoscopic character of $T$ associated to $H$*

$$\kappa_T : W(G,T) \backslash W(G(\mathbb{C}), T(\mathbb{C})) = H^1(\mathbb{R}, T) \longrightarrow \mu_2.$$

In what follows, we reformulate the above character identities to fit into our Main Theorem. For this purpose, we will use the explicit local Langlands correspondence for $Mp(2n)(\mathbb{R})$ à la J. Adams which we have discussed in Section 2 (see also [Ada10, Ada98]). Now we can start proving the expected endoscopic character identities for $(H, \tilde{G})$ at real places. We first transfer our Harish-Chandra parameter $\lambda'_p = (a_1, \cdots, a_{p_0}, -b_{q_0}, \cdots, -b_1)$ to

$$\lambda_n := (a_1, \cdots, a_{p_{p_0}}, b_{q_{q_0}}, \cdots, b_1)$$

using

$$w_0 = \prod_{j \in I_q} \sigma_j \in W^G_\mathbb{C}, \quad \text{with } I_q = \{q_1, \cdots, q_{q_0}\}.$$

Note that in Theorem 5.0.1, we have:

$$\begin{aligned}
Tran_{p_0,q_0}\Theta^H_\phi &= e^{\lambda'_p - \rho^\diamond}((1_{p_0}, -1_{q_0})) \sum_{w \in W^G_\mathbb{C}/W^K_\mathbb{C}} \kappa_T(w^{-1}) \Theta^{\tilde{G}}_{\tilde{\pi}_{w^{-1}\lambda'_p}} \\
&= e^{\lambda'_p - \rho^\diamond}((1_{p_0}, -1_{q_0})) \sum_{w \in W^G_\mathbb{C}/W^K_\mathbb{C}} \kappa_T(w^{-1}) \Theta^{\tilde{G}}_{\tilde{\pi}_{w^{-1}w_0\lambda_n}} \\
&= e^{\lambda'_p - \rho^\diamond}((1_{p_0}, -1_{q_0})) \sum_{w \in W^G_\mathbb{C}/W^K_\mathbb{C}} \kappa_T(w^{-1}w_0^{-1}) \Theta^{\tilde{G}}_{\tilde{\pi}_{w^{-1}\lambda_n}}.
\end{aligned} \quad (19)$$

Here $\rho^\diamond = (p_0, \cdots, 1) \times (q_0, \cdots, 1)$.

Note also that our expected character identity can be expressed as summation over $w \in W^G_\mathbb{C}/W^K_\mathbb{C}$ via the local Langlands correspondence for $\tilde{G}$ as follows.

$$\begin{aligned}
Tran_{p_0,q_0}\Theta^H_\phi &\stackrel{?}{=} \epsilon(1/2, \phi^{s=-1}, \psi) \sum_{\eta \in \hat{S}_\phi} \eta(s) \Theta^{\tilde{G}}_{\tilde{\pi}_\eta} \\
&= \epsilon(1/2, \phi^{s=-1}, \psi) \sum_{w \in W^G_\mathbb{C}/W^K_\mathbb{C}} \chi_{x_b^{-1}}(s) \chi_{w^{-1}}(s) \Theta^{\tilde{G}}_{\tilde{\pi}_{w^{-1}\lambda_n}}.
\end{aligned} \quad (20)$$

Here $s \in S_\phi$ with $s_{p_i} = 1, s_{q_j} = -1$, and $\chi_{w^{-1}}$ is the character of $S_\phi$ associated to $w^{-1}$ under the LLC for $\tilde{G}$ in Table 2.2.

Comparing the Formulas (19) and (20), it suffices to prove that

$$e^{\lambda'_p - \rho^\diamond}((1_{p_0}, -1_{q_0})) \kappa_T(w^{-1}w_0^{-1}) = \epsilon(1/2, \phi^{s=-1}, \psi) \chi_{x_b^{-1}}(s) \chi_{w^{-1}}(s) \quad (21)$$

holds for any $w \in W^G_\mathbb{C}/W^K_\mathbb{C}$, i.e.

$$e^{\lambda'_p - \rho^\diamond}((1_{p_0}, -1_{q_0})) \kappa_T(ww_0^{-1}) = \epsilon(1/2, \phi^{s=-1}, \psi) \chi_{x_b^{-1}}(s) \chi_w(s) \quad (22)$$

holds for any $w = \prod_{i \in I'_p} \sigma_i \times \prod_{j \in I'_q} \sigma_j \in W^G_\mathbb{C}$ with $I'_p \subset I_p$ and $I'_q \subset I_q$.

Observe that for $w = \prod_{i \in I'_p} \sigma_i \times \prod_{j \in I'_q} \sigma_j$ with $I'_p \subset I_p$ and $I'_q \subset I_q$, the corresponding character $\chi_\omega$ under LLC for $\tilde{G}$ has the following form:

$$\chi_{w,k} = (-1)^{n+1-k}, \quad k \notin I_p' \cup I_q';$$
$$\chi_{w,i} = (-1)^{n+2-i}, \quad i \in I_p';$$
$$\chi_{w,j} = (-1)^{n+2-j}, \quad j \in I_q'.$$

Here $\chi_{w,\cdot}$ stands for the corresponding component of $\chi_w$. Thus we obtain

$$\chi_w(s) = \prod_{j \in I_q \setminus I_q'} (-1)^{n+1-j} \prod_{j \in I_q'} (-1)^{n+2-j}$$
$$= (-1)^{\#I_q'} \prod_{j \in I_q} (-1)^{n+1-j}.$$

Also we have

$$\chi_{x_b^{-1}}(s) = \prod_{j \in I_q} (-1)^{n+1-j},$$

and by [Kna94, Theorem 3]), we have

$$\epsilon(1/2, \phi^{s=-1}, \psi) = \prod_{j=1}^{q_0} (-1)^{b_{q_j}+1/2}.$$

Therefore the right-hand side of Formula (22) equals

$$(-1)^{\#I_q'} \prod_{j=1}^{q_0} (-1)^{b_{q_j}+1/2}.$$

For the left-hand side, we know

$$w w_0^{-1} = \prod_{i \in I_p'} \sigma_i \times \prod_{j \in I_q \setminus I_q'} \sigma_j,$$

so

$$\kappa_T(w w_0^{-1}) = (-1)^{\#(I_q \setminus I_q')}.$$

Also observe that the $\widetilde{SL}_2(\mathbb{R})$-computation in [Ada98, P.175] gives

$$exp\left(\begin{pmatrix} & \pi \\ -\pi & \end{pmatrix}\right) = (-1_2, 1) \in \widetilde{SL}_2(\mathbb{R}), \qquad e^{-b_{q_j}}\left(exp\left(\begin{pmatrix} & \pi \\ -\pi & \end{pmatrix}\right)\right) = e^{-b_{q_j}\pi i} = (-1)(-1)^{b_{q_j}},$$

and the canonical element $-1_{2q_0} \in \widetilde{Sp}(2q_0)$ acting on Weil representation by -1 is

$$\begin{pmatrix} -1_{q_0} & \\ & -1_{q_0} \end{pmatrix} = \left(\begin{pmatrix} -1_{q_0} & \\ & -1_{q_0} \end{pmatrix}, \gamma((-1)^{q_0}, \psi)\right) \quad ([\text{Kud96, Chapter II Proposition 4.3}]).$$

So we have

$$e^{\lambda_p' - \rho^\diamond}((1_{p_0}, -1_{q_0})) = \prod_{j=1}^{q_0} (-1)(-1)^{b_{q_j}} \cdot (-1)^{\frac{q_0(q_0+1)}{2}} \gamma((-1)^{q_0}, \psi)$$
$$= \prod_{j=1}^{q_0} (-1)(-1)^{b_{q_j}} \cdot (-1)^{\frac{q_0(q_0+1)}{2}} (-1)^{\frac{q_0(q_0-1)}{2}} (-1)^{\frac{-q_0}{2}}$$
$$= \prod_{j=1}^{q_0} (-1)(-1)^{b_{q_j}} \cdot (-1)^{q_0} (-1)^{\frac{-q_0}{2}}$$
$$= \prod_{j=1}^{q_0} (-1)^{b_{q_j}-1/2}.$$

Thus the left-hand side of Formula (22) equals

$$(-1)^{\#(I_q \setminus I_q')} \prod_{j=1}^{q_0} (-1)^{b_{q_j}-1/2} = (-1)^{\#(I_q')} \prod_{j=1}^{q_0} (-1)^{b_{q_j}+1/2},$$

which really equals the right-hand side of Formula (22), whence the expected character identity holds.

**Corollary 5.0.2.** *Use the same notations defined above, we then have*

$$Tran_{p_0,q_0}(\Theta_\phi^H) = \epsilon(1/2, \phi^{s=-1}, \psi) \sum_{\eta \in \hat{S}_\phi} \eta(s) \Theta_{\tilde{\pi}_\eta}^{\tilde{G}}.$$

*Remark* 3. As noted, the case that $\lambda_p$ is not $G$-regular has not been discussed so far. But D. Renard has also established the general case using a coherent continuation argument in [Ren99]. Indeed, we could prove this via global argument as what we would do in the next section for non-archimedean cases. We leave the details for the reader.

*Remark* 4. As pointed out by the referee, the corollary follows easily from J. Adams 1994 ICM report [Ada95, Theorem 4.9] and W.W. Li's reformulation of Theorem 5.0.1 [Liar, Theorem 7.3.7]. But we would like to keep the above computation in its own right.

<u>Complex groups, i.e. $F = \mathbb{C}$.</u> Recall that $\tilde{G}(\mathbb{C}) = G(\mathbb{C}) \times \mu_8$, $H(\mathbb{C}) = SO_{2n'+1}(\mathbb{C}) \times SO_{2n''+1}(\mathbb{C})$. The endoscopy theory for $(\tilde{G}(\mathbb{C}), H(\mathbb{C}))$ shown by W.W. Li gives:

- The transfer factor $\Delta(\gamma, (t, \delta)) = t$ for all $t \in \mu_8$ and norm pair $(\gamma, \delta)$ with $\gamma \in H(\mathbb{C})$, $\delta \in G(\mathbb{C})$. This means that one may identify $C_c^\infty(\tilde{G}(\mathbb{C}))_{--}$ with $C_c^\infty(G(\mathbb{C}))$ via $\tilde{f} \longmapsto f(\cdot) := \tilde{f}(\cdot, 1)$.
- (Transfer map) For $f \in C_c^\infty(G(\mathbb{C}))$, there exists $f^H \in C_c^\infty(H(\mathbb{C}))$ such that
$$O_\gamma(f^H) = O_\delta(f)$$
for all norm pairs $(\gamma, \delta) \in H_{G-reg} \times G_{reg}$ with $\gamma = (\gamma', \gamma'') \in H_{G-reg}$.
- (Endoscopic character identity) For matching pairs $\big(\pi, (\sigma', \sigma'')\big)$ of representations of torus, i.e. $\pi = \sigma' \times \sigma''$, and matching pairs $(f, f^H)$ of test functions as above, we have

$$\epsilon(1/2, \sigma'') Trace\ Ind^G(\pi)(f) = \sigma''(-1) Trace\ Ind^G(\pi)(f) = Trace\ Ind^H((\sigma', \sigma''))(f^H),$$

which follows from the Weyl integration formula in the same way as in Proposition 4.1.2. Here the $(-1)$-factor in the norm correspondence results in a $\sigma''(-1)$-factor.

## 6. Transfer of discrete series characters: non-archimedean case

In this section, based on the endoscopic character identities for $\tilde{G}(\mathbb{R})$, we shall use stable trace formula of $\tilde{G}$ to deduce the desired endoscopic character identities for the tempered L-packets of $\tilde{G}$ which factors through endoscopic groups $H = H_{n',n''} := SO(2n'+1) \times SO(2n''+1)$ as discrete series L-packets. In what follows, we should discuss the cases $n' \neq n''$ and the case $n' = n''$ separately.

Before turning to the local-global argument, we first summarize some necessary inputs.

### 6.1. Globalization.
The input we need is the globalization of a local discrete L-parameter with prescribed local properties.

**Lemma 6.1.1** (Globalization of local field). *([Art13, Lemma 6.2.1]) Let $F$ be a $p$-adic field, and $r_0$ be a positive integer. Then there is a totally real number field $\dot{F}$ with the following properties:*

- $\dot{F}_{\nu_0} = F$ *for some finite place $\nu_0$ of $\dot{F}$.*
- *There are at least $r_0$ Archimedean places on $\dot{F}$.*

W may globalize our local data $(G_{\nu_0}, \psi_{\nu_0}, \langle \cdot, \cdot \rangle_{\nu_0})$ to be $(G, \psi, \langle \cdot, \cdot \rangle)$ as in [Liar, Proposition 8.4.1]. Let $\mathbb{H} = SO(2m+1)$ be a split special orthogonal group defined over $\mathbb{Z}$, and $S$ be a non-empty finite set of finite places of $\dot{F}$. Recall that $\dot{F}_S = \prod_{\nu \in S} \dot{F}_\nu$, and $\mathbb{H}(\dot{F}_S) = \prod_{\nu \in S} \mathbb{H}(\dot{F}_\nu)$. We write $\widehat{\mathbb{H}(\dot{F}_S)}$ as the unitary dual of $\mathbb{H}(\dot{F}_S)$ equipped with the Fell topology. We take $\hat{U}$ to be a $\hat{\mu}_S^{pl}$-regular relatively quasi-compact subset of $\widehat{\mathbb{H}(\dot{F}_S)}$ with $\hat{\mu}_S^{pl}(\hat{U}) > 0$, here $\hat{\mu}_S^{pl}$ is the Plancherel measure on $\widehat{\mathbb{H}(\dot{F}_S)}$.

**Theorem 6.1.2** (Globalization of local representations). *([Art13, Shi12]) Let $\mathbb{H}$, $S$ and $\hat{U}$ be as defined above. Then there exists infinitely many cuspidal automorphic representations $\pi$ of $\mathbb{H}(\mathbb{A}_{\dot{F}})$ such that*

- $\pi^{S,\infty}$ *is unramified,*
- $\pi_S \in \hat{U}$ *and*
- $\pi_\infty$ *are discrete series representations for all real places.*

*Proof.* One may refer to [Art13, Lemma 6.2.2] or [Shi12, Theorem 5.8] for details. □

In what follows, we take $S$ to be a sufficiently large finite subset of $V_{\dot{F}}$ consisting of finite places of $\dot{F}$. Fix a discrete series representation $\pi_{\nu_i}$ for each $\nu_i \in S$, and fix a discrete series L-parameter $\phi_{\infty_i}$ in general position in the sense of Arthur for each $\infty_i \in V_\infty$ (see [Art13, Lemma 6.2.2]).

**Corollary 6.1.3** (Globalization of L-parameters). *There exists an elliptic tempered global A-parameter $\Psi$ of $\mathbb{H}$ such that*
- $\Psi_{\nu_i}$ *is the discrete series L-parameter associated to $\pi_{\nu_i}$ for each $\nu_i \in S$,*
- $\Psi_{\infty_i}$ *is the fixed discrete series L-parameter $\phi_{\infty_i}$ for each $\infty_i \in V_\infty$ and*
- $\Psi$ *is unramified at the remaining finite places.*

*Proof.* This follows from [Art13, Corollary 6.2.3]. □

Now we may proceed to the proof of the Main Theorem. Note that our spherical fundamental lemma for $\tilde{G}$ is not available at 2-adic places, thus the proof shall be further separated into two cases, i.e. $\nu_0 \nmid 2$ and $\nu_0 \mid 2$.

6.2. **Cases $n' \neq n''$ and $\nu_0 \nmid 2$.** Our main strategy is to apply twice the simple stable trace formula II of $\tilde{G}$, i.e. Corollary 3.2.5, to get an auxiliary equation and a target one respectively. Our Main Theorem would then be derived from an easy comparison of these two equations.

**Picking test functions I: $\nu_0 \notin S$ and $\nu_0 \nmid 2$.** Fix a sufficiently large $S$ including all ramified places of $(G,\psi)$ and places over 2 except $\nu_0$. Let $\infty_1$ be a fixed archimedean place of $\dot{F}$. To apply Gan–Ichino's multiplicity formula, one should first construct an elliptic tempered A-parameter of $\tilde{G}$ such that its local component at $\nu_0$ is the same as in our Main Theorem. To achieve this, one may apply Corollary 6.1.3 to $H_1 = SO(2n'+1)$ and $H_2 = SO(2n''+1)$ separately, but with different local constraints as follows.

- At a finite place $\nu_1 \in S$, we fix an irreducible discrete series representation $\pi_1$ of $H_1$, and an irreducible discrete series representation $\pi_2$ of $H_2$, such that both are associated to simple L-parameters.
- At each real place $\infty_i \in V_\infty$, we also fix a discrete series L-parameter of $H_1$ and $H_2$ respectively, and require their composition as L-parameter of $\tilde{G}$ is $G$-regular (see [Liar, Section 7.3]).

As $n' \neq n''$, by Corollary 6.1.3, we obtain an elliptic tempered A-parameter $\Psi$ of $\tilde{G}$ which factors through $H_1 \times H_2$ as $\Psi_1 \times \Psi_2$ in a unique way up to $\hat{H}_1 \times \hat{H}_2$-conjugation.

Now turn to the choice of test functions. In order to apply the simple stable trace formula, we shall further specify the test functions we will be using. We consider $\tilde{f} = \otimes_\nu \tilde{f}_\nu \in C^\infty_{c,--}(\tilde{G}(\mathbb{A}_{\dot{F}}))$ satisfying

- for finite places $\nu \notin S$, $\tilde{f}_\nu$ belongs to the spherical Hecke algebra, and $f^H_\nu$ is the corresponding spherical function of $H$. This is guaranteed by the spherical fundamental lemma for $(\tilde{G}, H)$ (cf. [Luoar]).
- for the finite place $\nu_1 \in S$, we take a cuspidal test function $\tilde{f}_{\nu_1}$ such that it transfers to $f^H_{\nu_1}$ with $\Theta^H_{\Psi_{\nu_1}}(f^H_{\nu_1}) \neq 0$.
- for another finite place $\nu_2 \in S$, we take a cuspidal test function $\tilde{f}_{\nu_2}$ such that it transfers to $f^H_{\nu_2}$ for the fixed endoscopic group $H$, and transfers to 0 for other elliptic endoscopic groups. We may further require that $\Theta^H_{\Psi_{\nu_2}}(f^H_{\nu_2}) \neq 0$ by the choice of an appropriate discrete series representation $\pi_{\nu_2}$.
- for an extra finite place $\nu_3 \in S$, we take a cuspidal test function $\tilde{f}_{\nu_3}$ such that it is supported on the regular semisimple elliptic locus of $\tilde{G}_{\nu_3}$, and such that it transfers to $f^H_{\nu_3}$ with $\Theta^H_{\Psi_{\nu_3}}(f^H_{\nu_3}) \neq 0$. This is granted by Hales' Lemma in [Hal95, Lemma 5.1] and Arthur's stable density theorem (see [Art13, Proposition 6.5.1], and refer to [Sha90, Conjecture 9.2] for the explicit statement) as follows: From Hales' construction in [Hal95, Lemma 5.1], there exist cuspidal functions supported on the regular semisimple elliptic locus of $\tilde{G}_{\nu_3}$ such that their $\kappa$-orbital integrals are ordinary orbital integrals ([Kot86, Proposition 7.1]). For such a test function of $\tilde{G}$, by Theorem 3.1.1(iii), we then get a cuspidal test function on $H$ which is non-zero as a stable distribution. Further by Arthur's stable density theorem [Art13, Proposition 6.5.1], we know that $\Theta^H_{\Psi_{\nu_3}}(f^H_{\nu_3}) \neq 0$ for some discrete L-packet.

- for other $\nu \in S$, we take $(\tilde{f}_\nu, f^H_\nu)$ to be a transfer pairing for $(\tilde{G}, H)$ such that $\Theta^H_{\Psi_\nu}(f^H_\nu) \neq 0$.
- for $\infty_1$, set
$$\tilde{f}_{\infty_1} = \frac{1}{|S_{\Psi_{\infty_1}}|} \sum_{\eta \in \hat{S}_{\Psi_{\infty_1}}} \eta(s)\tilde{f}_\eta$$
with $\tilde{f}_\eta$ pseudo-coefficients of the representations $\tilde{\pi}_\eta$ in the L-packet of $\Psi_{\infty_1}$.
- for other archimedean places $\infty_i$, we take $\tilde{f}_{\infty_i}$ to be a pseudo-coefficient of $\tilde{\pi}_{\eta_i}$ for a fixed $\eta_i \in \hat{S}_{\Psi_{\infty_i}}$.

For simplicity, we denote $A_{\Psi_S} := \prod_{\nu \in S} A_{\Psi_\nu}$ to be the product of component groups of $\psi_\nu$, and

$$\epsilon(1/2, \Psi_S^{s=-1}, \psi) = \prod_{\nu \in S} \epsilon(1/2, \Psi_\nu^{s=-1}, \psi), \qquad \eta_S = \prod_{\nu \in S} \eta_\nu \quad \text{and} \quad \Theta_{\eta_S}(\tilde{f}_S) = \prod_{\nu \in S} \Theta_{\eta_\nu}(\tilde{f}_\nu).$$

Applying the simple stable trace formula i.e. Corollary 3.2.5, we have:

$$\prod_{\infty_i \in V_\infty \setminus \{\infty_1\}} \Theta_{\eta_i}(\tilde{f}_{\infty_i}) \times \sum_{\substack{\eta_1, \eta_S: \\ \Delta^*(\prod_i \eta_i \prod_\nu \eta_\nu) = \epsilon_\Psi}} \Theta_{\eta_1}(\tilde{f}_{\infty_1}) \Theta_{\eta_S}(\tilde{f}_S)$$
(23)
$$= |S_\Psi|^{-1} \prod_{\infty_i \in V_\infty \setminus \{\infty_1\}} \Theta_{\Psi_{\infty_i}}(f^H_{\infty_i}) \times \Theta_{\Psi_{\infty_1}}(f^H_{\infty_1}) \prod_{\nu \in S} \Theta_{\Psi_\nu}(f^H_\nu).$$

Note that
$$\Theta_{\eta_1}(\tilde{f}_{\infty_1}) = \frac{1}{|S_{\Psi_{\infty_1}}|}\eta_1(s), \quad \Theta_{\eta_i}(\tilde{f}_{\infty_i}) = 1.$$

So the left-hand side of Formula (23) equals

$$\sum_{\substack{\eta_1, \eta_S: \\ \Delta^*(\prod_i \eta_i \prod_\nu \eta_\nu) = \epsilon_\Psi}} \Theta_{\eta_1}(\tilde{f}_{\infty_1}) \Theta_{\eta_S}(\tilde{f}_S)$$
$$= \sum_{\eta_S} \sum_{\substack{\eta_1: \\ \Delta^*(\prod_i \eta_i \prod_\nu \eta_\nu) = \epsilon_\Psi}} \Theta_{\eta_1}(\tilde{f}_{\infty_1}) \Theta_{\eta_S}(\tilde{f}_S)$$
$$= \sum_{\eta_S} \sum_{\substack{\eta_1: \\ \Delta^*(\prod_i \eta_i \prod_\nu \eta_\nu) = \epsilon_\Psi}} \frac{1}{|S_{\Psi_{\infty_1}}|} \eta_1(s) \Theta_{\eta_S}(\tilde{f}_S)$$
$$= \frac{|S_{\Psi_{\infty_1}}|}{|S_\Psi|} \frac{\epsilon(1/2, \Psi^{s=-1}, \psi)}{|S_{\Psi_{\infty_1}}|} \prod_{\substack{i: \\ \infty_i \in V_\infty \setminus \{\infty_1\}}} \eta_i(s) \sum_{\eta_S \in \hat{A}_{\Psi_S}} \eta_S(s) \Theta_{\eta_S}(\tilde{f}_S).$$

Here the ratio $\frac{|S_{\Psi_{\infty_1}}|}{|S_\Psi|}$ in the last equation results from the relation $\epsilon_\Psi : S_\Psi \hookrightarrow S_{\Psi_{\infty_1}} \xrightarrow{\eta_1} \mu_2$ and the fact that $S_{\Psi_{\infty_1}}$ is a finite abelian 2-group.

For the RHS of Formula 23, note that for $\infty_i \in V_\infty \setminus \{\infty_1\}$, by Corollary 5.0.2, we have
$$\Theta_{\Psi_{\infty_i}}(f^H_{\infty_i}) = \epsilon(1/2, \Psi_{\infty_i}^{s=-1}, \psi)\eta_i(s), \quad \Theta_{\Psi_{\infty_1}}(f^H_{\infty_1}) = \epsilon(1/2, \Psi_{\infty_1}^{s=-1}, \psi).$$

Thus the right-hand side of Formula (23) equals
$$|S_\Psi|^{-1} \prod_{\infty_i \in V_\infty \setminus \{\infty_1\}} \epsilon(1/2, \Psi_{\infty_i}^{s=-1}, \psi)\eta_i(s) \times \epsilon(1/2, \Psi_{\infty_1}^{s=-1}, \psi) \prod_{\nu \in S} \Theta_{\Psi_\nu}(f^H_\nu).$$

Combining both sides together, we have:

(**AUXILIARY**) $\qquad \epsilon(1/2, \Psi_S^{s=-1}, \psi) \sum_{\eta_S \in \hat{A}_{\Psi_S}} \eta_S(s) \Theta_{\eta_S}(\tilde{f}_S) = \prod_{\nu \in S} \Theta_{\Psi_\nu}(f^H_\nu).$

Note that this formula is meaningful as the right hand side is non-zero by the choice of test functions. This is our auxiliary equation which will be used later on.

**Picking test functions II:** $\nu_0 \notin S$ **and** $\nu_0 \nmid 2$ . Fix a new $S_{new} := \{\nu_0\} \cup S$ with $S$ being the old $S$ as above. For all $\nu \neq \nu_0$, we put the same restrictions as above, and impose the same conditions on the test functions $\tilde{f}_\nu$. To get a global A-parameter as above, we take the $\nu_0$-component to be our target discrete series L-parameter, and denote by $\Psi'$ the derived global A-parameter of $\tilde{G}$.

Run the same procedure for the new $S$ and $\Psi'$ with $\tilde{f}_{\nu_0}$ arbitrary, and write out terms associated to $S_{new}$ in terms of $\nu_0$ and the old $S$, we get an analogous formula as Formula (**AUXILIARY**), i.e. our target equation:

(**TARGET**)
$$\epsilon(1/2, \Psi'^{s=-1}_S, \psi)\epsilon(1/2, \Psi'^{s=-1}_{\nu_0}, \psi) \sum_{\eta_{\nu_0} \in \hat{A}_{\Psi'_{\nu_0}}} \eta_{\nu_0}(s)\Theta_{\eta_{\nu_0}}(\tilde{f}_{\nu_0}) \sum_{\eta_S \in \hat{A}_{\Psi_S}} \eta_S(s)\Theta_{\eta_S}(\tilde{f}_S)$$
$$= \Theta_{\Psi'_{\nu_0}}(f^H_{\nu_0}) \prod_{\nu \in S} \Theta_{\Psi'_\nu}(f^H_\nu).$$

Observe that
$$\epsilon(1/2, \Psi'^{s=-1}_S, \psi) = \epsilon(1/2, \Psi^{s=-1}_S, \psi).$$

By the non-vanishing property of Formula (**AUXILIARY**), we get, via an easy comparison of Formula (**AUXILIARY**) and Formula (**TARGET**),

$$\epsilon(1/2, \Psi'^{s=-1}_{\nu_0}, \psi) \sum_{\eta_{\nu_0} \in \hat{A}_{\Psi'_{\nu_0}}} \eta_{\nu_0}(s)\Theta_{\eta_{\nu_0}}(f^H_{\nu_0}) = \Theta_{\Psi'_{\nu_0}}(f^H_{\nu_0}).$$

Whence we have finished the proof of the Main Theorem for the cases $n' \neq n''$ and $\nu_0 \nmid 2$.

6.3. **Cases $n' = n''$ and $\nu_0 \nmid 2$.** Now let us come back to the case $n' = n''$ and $\nu_0 \nmid 2$. In this case, we write $H = H_1 \times H_2$ with $H_1 = H_2 = SO(n+1)$ (n is even), and denote our target discrete series L-parameter $\phi = \phi_1 \times \phi_2$. Note that if $\phi_1 \simeq \phi_2$, then we can adapt the same argument as in cases $n' \neq n''$ with minor modifications as follows.

- All the choices of test functions and fixed local L-parameters are the same as in cases $n' \neq n''$. Instead of getting global A-parameters $\Psi$ and $\Psi'$ factoring uniquely through $\hat{H}$, we find that $\Psi$ and $\Psi'$ are factoring through $\hat{H}$ in two ways, i.e. $\Psi = \Psi_1 \oplus \Psi_2$ and $\Psi^{op} = \Psi_2 \oplus \Psi_1$, similarly $\Psi' = \Psi'_1 \oplus \Psi'_2$ and $\Psi'^{op} = \Psi'_2 \oplus \Psi'_1$.
- Applying the simple stable trace formulas as above, instead of getting a equation involving only one term on the right hand side of Formula (23), we get

(24)
$$\prod_{\infty_i \in V_\infty \setminus \{\infty_1\}} \Theta_{\eta_i}(\tilde{f}_{\infty_i}) \times \sum_{\substack{\eta_1, \eta_S: \\ \Delta^*(\prod_i \eta_i \prod_\nu \eta_\nu) = \epsilon_\Psi}} \Theta_{\eta_1}(\tilde{f}_{\infty_1})\Theta_{\eta_S}(\tilde{f}_S)$$
$$= |S_\Psi|^{-1} \prod_{\infty_i \in V_\infty \setminus \{\infty_1\}} \Theta_{\Psi_{\infty_i}}(f^H_{\infty_i}) \times \Theta_{\Psi_{\infty_1}}(f^H_{\infty_1}) \prod_{\nu \in S} \Theta_{\Psi_\nu}(f^H_\nu)$$
$$+ |S_\Psi|^{-1} \prod_{\infty_i \in V_\infty \setminus \{\infty_1\}} \Theta_{\Psi^{op}_{\infty_i}}(f^H_{\infty_i}) \times \Theta_{\Psi^{op}_{\infty_1}}(f^H_{\infty_1}) \prod_{\nu \in S} \Theta_{\Psi^{op}_\nu}(f^H_\nu).$$

- Therefore our auxiliary equation should be

(25)
$$\epsilon(1/2, \Psi^{s=-1}, \psi) \prod_{\substack{i: \\ \infty_i \in V_\infty \setminus \{\infty_1\}}} \eta_i(s) \sum_{\eta_S \in \hat{A}_{\Psi_S}} \eta_S(s)\Theta_{\eta_S}(\tilde{f}_S)$$
$$= \prod_{\infty_i \in V_\infty \setminus \{\infty_1\}} \epsilon(1/2, \Psi^{s=-1}_{\infty_i}, \psi)\eta_i(s) \times \epsilon(1/2, \Psi^{s=-1}_{\infty_1}, \psi) \prod_{\nu \in S} \Theta_{\Psi_\nu}(f^H_\nu)$$
$$+ \prod_{\infty_i \in V_\infty \setminus \{\infty_1\}} \epsilon(1/2, \Psi^{op,s=-1}_{\infty_i}, \psi)\eta_i(-s) \times \epsilon(1/2, \Psi^{op,s=-1}_{\infty_1}, \psi) \prod_{\nu \in S} \Theta_{\Psi^{op}_\nu}(f^H_\nu).$$

- Similarly our target equation should be

$$\epsilon(1/2, \Psi'^{s=-1}, \psi) \prod_{\infty_i \in V_\infty \setminus \{\infty_1\}} \eta_i(s) \sum_{\eta_S \in \hat{A}_{\Psi_S}} \eta_S(s) \Theta_{\eta_S}(\tilde{f}_S) \times \sum_{\eta_{\nu_0}} \eta_{\nu_0} \Theta_{\eta_{\nu_0}}(\tilde{f}_{\nu_0})$$

$$(26) \quad = \prod_{\infty_i \in V_\infty \setminus \{\infty_1\}} \epsilon(1/2, \Psi'^{s=-1}_{\infty_i}, \psi)\eta_i(s) \times \epsilon(1/2, \Psi'^{s=-1}_{\infty_1}, \psi) \prod_{\nu \in S} \Theta_{\Psi'_\nu}(f^H_\nu) \times \Theta_{\Psi'_{\nu_0}}(f^H_{\nu_0})$$

$$+ \prod_{\infty_i \in V_\infty \setminus \{\infty_1\}} \epsilon(1/2, \Psi'^{op,s=-1}_{\infty_i}, \psi)\eta_i(-s) \times \epsilon(1/2, \Psi'^{op,s=-1}_{\infty_1}, \psi) \prod_{\nu \in S} \Theta_{\Psi'^{op}_\nu}(f^H_\nu) \times \Theta_{\Psi'^{op}_{\nu_0}}(f^H_{\nu_0}).$$

- Note that $\Theta_{\Psi'_{\nu_0}}(f^H_{\nu_0}) = \Theta_{\Psi'^{op}_{\nu_0}}(f^H_{\nu_0})$ as $\phi_1 \simeq \phi_2$, and

$$\epsilon(1/2, \Psi^{s=-1}, \psi)\epsilon(1/2, \phi_2, \psi_{\nu_0}) = \epsilon(1/2, \Psi'^{s=-1}, \psi).$$

Thus by an easy comparison, we obtain

$$\epsilon(1/2, \Psi'^{s=-1}_{\nu_0}, \psi) \sum_{\eta_{\nu_0} \in \hat{A}_{\Psi'_{\nu_0}}} \eta_{\nu_0}(s) \Theta_{\eta_{\nu_0}}(f^H_{\nu_0}) = \Theta_{\Psi'_{\nu_0}}(f^H_{\nu_0}),$$

which is exactly what the Main Theorem says.

Note that only the cases $n' = n''$ and $\phi = \phi_1 \oplus \phi_2$ with $\phi_1 \not\simeq \phi_2$ have not been touched. Fortunately, the remaining cases also follow from the same argument adopted in the cases $n' \neq n''$, but with another minor modification as follows:

- All the choices of test functions are the same as in cases $n' \neq n''$ except at a finite place $\nu_4 \in S$. At such a place, we fix the local $L$-parameter of $\tilde{G}$ as $\phi_{\nu_4} = \phi'_{\nu_4} \oplus \phi''_{\nu_4}$ with simple $L$-parameters $\phi'_{\nu_4} \not\simeq \phi''_{\nu_4}$, and choose a test function $f^H_{\nu_4} = f^{H_1}_{\nu_4} \times f^{H_2}_{\nu_4}$, such that $f^{H_1}_{\nu_4}$ is the pseudo-coefficient of the $L$-parameter $\phi'_{\nu_4}$. Then by Theorem 3.1.1(iii), we take a test function $\tilde{f}_{\nu_4}$ of $\tilde{G}$ such that $(\tilde{f}_{\nu_4}, f^H_{\nu_4})$ is a transfer pair for $(\tilde{G}, H)$. We also find that $\Psi$ and $\Psi'$ are factoring through $\mathbb{H}$ in two ways, i.e. $\Psi = \Psi_1 \oplus \Psi_2$ and $\Psi^{op} = \Psi_2 \oplus \Psi_1$, similarly $\Psi' = \Psi'_1 \oplus \Psi'_2$ and $\Psi'^{op} = \Psi'_2 \oplus \Psi'_1$.
- Applying the simple stable trace formulas as above, instead of getting an equation involving two terms on the right hand side of Formula (24), we get only one term as in the cases $n' \neq n''$:

$$\prod_{\infty_i \in V_\infty \setminus \{\infty_1\}} \Theta_{\eta_i}(\tilde{f}_{\infty_i}) \times \sum_{\substack{\eta_1, \eta_S: \\ \Delta^*(\prod_i \eta_i \prod_\nu \eta_\nu) = \epsilon_\Psi}} \Theta_{\eta_1}(\tilde{f}_{\infty_1}) \Theta_{\eta_S}(\tilde{f}_S)$$

$$= |S_\Psi|^{-1} \prod_{\infty_i \in V_\infty \setminus \{\infty_1\}} \Theta_{\Psi_{\infty_i}}(f^H_{\infty_i}) \times \Theta_{\Psi_{\infty_1}}(f^H_{\infty_1}) \prod_{\nu \in S} \Theta_{\Psi_\nu}(f^H_\nu).$$

This results from the vanishing condition:

$$\Theta_{\Psi^{op}_{\nu_4}}(f^H_{\nu_4}) = \Theta^{H_1}_{\phi''_{\nu_4}}(f^{H_1}_{\nu_4}) \Theta^{H_2}_{\phi'_{\nu_4}}(f^{H_2}_{\nu_4}) = 0.$$

- Thus we get an auxiliary equation and a target equation as follows.

$$\epsilon(1/2, \Psi^{s=-1}_S, \psi) \sum_{\eta_S \in \hat{A}_{\Psi_S}} \eta_S(s) \Theta_{\eta_S}(\tilde{f}_S) = \prod_{\nu \in S} \Theta_{\Psi_\nu}(f^H_\nu), \text{ and}$$

$$(27) \quad \epsilon(1/2, \Psi'^{s=-1}_S, \psi)\epsilon(1/2, \Psi'^{s=-1}_{\nu_0}, \psi) \sum_{\eta_{\nu_0} \in \hat{A}_{\Psi'_{\nu_0}}} \eta_{\nu_0}(s) \Theta_{\eta_{\nu_0}}(f^H_{\nu_0}) \sum_{\eta_S \in \hat{A}_{\Psi_S}} \eta_S(s) \Theta_{\eta_S}(\tilde{f}_S)$$

$$= \Theta_{\Psi'_{\nu_0}}(f^H_{\nu_0}) \prod_{\nu \in S} \Theta_{\Psi'_\nu}(f^H_\nu).$$

Therefore we obtain

$$\epsilon(1/2, \Psi'^{s=-1}_{\nu_0}, \psi) \sum_{\eta_{\nu_0} \in \hat{S}_{\Psi'_{\nu_0}}} \eta_{\nu_0}(s) \Theta_{\eta_{\nu_0}}(f^H_{\nu_0}) = \Theta_{\Psi'_{\nu_0}}(f^H_{\nu_0}).$$

So far, we have finished the proof of the Main Theorem for all cases except $\nu_0 \mid 2$.

*Remark* **5.** *In view of the above arguments, we realize that it is not necessary to prove the Main Conjecture for the cases $n' \neq n''$ and $n' = n''$ separately. But we keep it there as the arguments are interesting in their own right.*

6.4. **Cases $\nu_0 \mid 2$.** In this case, we require further that our global field $\dot{F}/\mathbb{Q}$ with $\dot{F}_{\nu_0} = F$, such that $\nu_0$ is the only prime over 2. Given the endoscopic character identities for all places except $\nu_0$, we may use only the target equation to deduce the Main Theorem for $\nu_0$. Following the same procedure as in the case $n' \neq n''$, we pick up a large enough set $S$ of finite places of $\dot{F}$ including all ramified places except $\nu_0$, and construct an elliptic tempered A-parameter $\Psi$ of $H$ as before. Another data is about test functions. We take the same test functions $\tilde{f} = \otimes_\nu \tilde{f}_\nu \in C_{c,--}^\infty(\tilde{G}(\mathbb{A}_{\dot{F}}))$ as in the last case except $\tilde{f}_{\nu_0}$ being arbitrary. Routinely, we also get the Formula (**TARGET**).

$$\epsilon(1/2, \Psi'^{s=-1}_S, \psi)\epsilon(1/2, \Psi'^{s=-1}_{\nu_0}, \psi) \sum_{\eta_{\nu_0} \in \hat{A}_{\Psi'_{\nu_0}}} \eta_{\nu_0}(s)\Theta_{\eta_{\nu_0}}(\tilde{f}_{\nu_0}) \sum_{\eta_S \in \hat{A}_{\Psi_S}} \eta_S(s)\Theta_{\eta_S}(\tilde{f}_S)$$
$$= \Theta_{\Psi'_{\nu_0}}(f^H_{\nu_0}) \prod_{\nu \in S} \Theta_{\Psi'_\nu}(f^H_\nu).$$

Thus we may readily conclude the endoscopic character identity for $\nu_0$ as the canceled terms are non-zero. Whence we have completely finished the proof of the Main Theorem for all cases.

## References

[AB98] Jeffrey Adams and Dan Barbasch, *Genuine representations of the metaplectic group*, Compositio Mathematica **113** (1998), no. 01, 23–66.

[ABV12] Jeffrey Adams, Dan Barbasch, and David A Vogan, *The Langlands classification and irreducible characters for real reductive groups*, vol. 104, Springer Science & Business Media, 2012.

[Ada95] Jeffrey Adams, *Genuine representations of the metaplectic group and epsilon factors*, Proceedings of the International Congress of Mathematicians, Springer, 1995, pp. 721–731.

[Ada98] ______, *Lifting of characters on orthogonal and metaplectic groups*, Duke Mathematical Journal **92** (1998), no. 1, 129–178.

[Ada10] ______, *Discrete series and characters of the component group*, http://www.math.umd.edu/~jda/preprints/discreteSeriesSigns.total.pdf, 2010.

[Art88] James Arthur, *The invariant trace formula. II. Global theory*, Journal of the American Mathematical Society **1** (1988), no. 3, 501–554.

[Art05] ______, *An introduction to the trace formula*, Harmonic analysis, the trace formula, and Shimura varieties **4** (2005), 1–263.

[Art13] ______, *The Endoscopic Classification of Representations Orthogonal and Symplectic Groups*, vol. 61, American Mathematical Soc., 2013.

[Bor79] Armand Borel, *Automorphic l-functions*, Automorphic forms, representations and L-functions (Proc. Sympos. Pure Math., Oregon State Univ., Corvallis, Ore., 1977), Part, vol. 2, 1979, pp. 27–61.

[CG15] Ping-Shun Chan and Wee Teck Gan, *The local Langlands conjecture for GSp (4) III: Stability and twisted endoscopy*, Journal of Number Theory **146** (2015), 69–133.

[CG16] Kwangho Choiy and David Goldberg, *Invariance of R-groups between p-adic inner forms of quasi-split classical groups*, Transactions of the American Mathematical Society **368** (2016), no. 2, 1387–1410.

[GGP12] Wee Teck Gan, Benedict H Gross, and Dipendra Prasad, *Symplectic local root numbers, central critical L-values, and restriction problems in the representation theory of classical groups*, Astérisque **346** (2012), 1–109.

[GI17] Wee Teck Gan and Atsushi Ichino, *The Shimura-Waldspurger correspondence for $Mp_{2n}$*, arXiv preprint arXiv:1705.10106 (2017).

[GL16] Wee Teck Gan and Wen-Wei Li, *The Shimura-Waldspurger correspondence for $Mp(2n)$*, Preprint (2016).

[Gol94] David Goldberg, *Reducibility of induced representations for Sp(2n) and SO(n)*, American Journal of Mathematics **116** (1994), no. 5, 1101–1151.

[GP92] Benedict H Gross and Dipendra Prasad, *On the decomposition of a representation of $SO_n$ when restricted to $SO_{n-1}$*, Canad. J. Math **44** (1992), 974–1002.

[GS12] Wee Teck Gan and Gordan Savin, *Representations of metaplectic groups I: epsilon dichotomy and local Langlands correspondence*, Compositio Mathematica **148** (2012), no. 06, 1655–1694.

[Hal95] Thomas C Hales, *On the fundamental lemma for standard endoscopy: reduction to unit elements*, Canadian Journal of Mathematics **47** (1995), no. 5, 974–994.

[How10] Tatiana K Howard, *Lifting of characters on p-adic orthogonal and metaplectic groups*, Compositio Mathematica **146** (2010), no. 03, 795–810.

[Kna94] A. W. Knapp, *The local Langlands correspondence: the non-Archimedean case*, Proc. Symp. Pure Math, vol. 55, Motives, 1994, pp. 393–410.

[Kot86] Robert E Kottwitz, *Stable trace formula: elliptic singular terms*, Mathematische Annalen **275** (1986), no. 3, 365–399.
[Kud96] Stephen Kudla, *Notes on the local theta correspondence*, http://www.math.toronto.edu/skudla/castle.pdf, 1996.
[Li11] Wen-Wei Li, *Transfert des intégrales orbitales pour le groupe métaplectique*, Compositio Mathematica **147** (2011), no. 02, 524–590.
[Li12a] ______, *La formule des traces pour les revêtements de groupes réductifs connexes. II. Analyse harmonique locale*, Annales scientifiques de l'ENS **45** (2012), no. 05.
[Li12b] ______, *Le lemme fondamental pondéré pour le groupe métaplectique*, Canadian Journal of Mathematics **64** (2012), no. 03, 497–543.
[Li13] ______, *La formule des traces pour les revêtements de groupes réductifs connexes III: Le développement spectral fin*, Mathematische Annalen **356** (2013), no. 3, 1029–1064.
[Li14a] ______, *La formule des traces pour les revêtements de groupes réductifs connexes. I. Le développement géométrique fin*, Journal für die reine und angewandte Mathematik (Crelles Journal) **2014** (2014), no. 686, 37–109.
[Li14b] ______, *La formule des traces pour les revêtements de groupes réductifs connexes. IV. Distributions invariantes*, Annales de l'Institut Fourier **64** (2014), no. 06, 2379–2448.
[Li15] ______, *La formule des traces stable pour le groupe métaplectique: les termes elliptiques*, Inventiones Mathematicae **202** (2015), no. 02, 743–838.
[Liar] ______, *Spectral transfer for metaplectic groups. I. Local character relations*, Journal of the Institute of Mathematics of Jussieu (To appear).
[Luoar] Caihua Luo, *Spherical fundamental lemma for metaplectic groups*, Canadian Journal of Mathematics (To appear).
[MR17] Colette Moeglin and David Renard, *Sur les paquets dArthur des groupes classiques et unitaires non quasi-déployés*, 2017.
[Ren99] David Renard, *Endoscopy for Mp(2n, $\mathbb{R}$)*, American Journal of Mathematics **121** (1999), no. 6, 1215–1244.
[Sch98] Jason Paul Schultz, *Lifting of Characters of $\widetilde{SL}_2(F)$ and $SO_{1,2}(F)$ for F A Nonarchimedean Local Field*, Ph.D. thesis, University of Maryland, College Park, Md., 1998.
[Sha90] Freydoon Shahidi, *A proof of Langlands' conjecture on Plancherel measures: complementary series of p-adic groups*, Annals of Mathematics **132** (1990), no. 2, 273–330.
[Shi12] Sug Woo Shin, *Automorphic Plancherel density theorem*, Israel Journal of Mathematics **192** (2012), no. 1, 83–120.
[Wal80] Jean Loup Waldspurger, *Correspondance de Shimura*, J. Math. Pures Appl **59** (1980), no. 1, 1–133.
[Wal91] Jean-Loup Waldspurger, *Correspondances de Shimura et quaternions*, Forum Mathematicum, vol. 3, 1991, pp. 219–308.
[Wal01] ______, *Intégrales orbitales nilpotentes et endoscopie pour les groupes classiques non ramifiés*, Astérisque **269** (2001).

DEPARTMENT OF MATHEMATICS, NATIONAL UNIVERSITY OF SINGAPORE, 10 LOWER KENT RIDGE ROAD SINGAPORE 119076
*E-mail address*: cluo@u.nus.edu